\documentclass[12pt,a4paper]{amsproc}
\usepackage{a4}
\usepackage{graphicx} 
\usepackage{amssymb}
\usepackage{amscd}
\usepackage[mathscr]{eucal}
\usepackage{amsmath,tikz-cd}
\usepackage{amsthm}
\usepackage{mathtools}
\usepackage{hyperref}
\usepackage{comment}
\usepackage[text={6in,8.5in},centering]{geometry}

\newtheorem{thm}{Theorem}[section]
\newtheorem{cor}[thm]{Corollary}
\newtheorem{lemma}[thm]{Lemma}
\newtheorem{prop}[thm]{Proposition}

\newtheorem{claim}{Claim}

\newcommand{\Id}{\operatorname{Id}}

\theoremstyle{definition}
\newtheorem{defi}[thm]{Definition}
\newtheorem{remark}[thm]{Remark}

\newcommand{\pB}{\operatorname{\textbf{pB}}}
\newcommand{\QB}{\operatorname{\textbf{QB}}}
\newcommand{\pBL}{\operatorname{\textbf{pBL}}}
\newcommand{\pBLp}{\operatorname{\textbf{pBL}}^{\textbf{(p)}}}
\newcommand{\dens}{\operatorname{dens}}
\newcommand{\Hom}{\operatorname{Hom}}
\newcommand{\Fin}{\operatorname{Fin}}
\newcommand{\FVL}{\operatorname{FVL}}
\newcommand{\FBL}{\operatorname{FBL}}
\newcommand{\FpBL}{\operatorname{FpBL}}
\newcommand{\FxBL}[1]{\operatorname{F{#1}BL}^{(#1)}}
\newcommand{\Fx}[1]{\operatorname{F{#1}BL}}
\newcommand{\FqL}{\operatorname{FQBL^{(L)}}}
\newcommand{\FBLp}{\operatorname{FBL^{(p)}}}
\newcommand{\Fq}{\operatorname{FQBL}}

\renewcommand{\Re}{\operatorname{Re}}

\newcommand{\FpBLp}{\operatorname{FpBL}^{(p)}}
\newcommand{\norm}[1]{{\left\vert\kern-0.25ex\left\vert #1 
		\right\vert\kern-0.25ex\right\vert}}
\newcommand{\nnorm}[1]{{\left\vert\kern-0.25ex\left\vert\kern-0.25ex\left\vert #1 
		\right\vert\kern-0.25ex\right\vert\kern-0.25ex\right\vert}}
\newcommand{\R}{\mathbb R}
\newcommand\restr[2]{{
		\left.\kern-\nulldelimiterspace 
		#1 
		\right|_{#2} 
}}

\title[Free quasi-Banach lattices]{Free quasi-Banach lattices}

\author[Salguero-Alarc\'on]{Alberto Salguero-Alarc\'on}
\address{Departamento de An\'alisis Matem\'atico y Matem\'atica Aplicada\\
Universidad Complutense de Madrid\\
Plaza de las Ciencias 3\\
28040 Madrid, Spain.}
\email{albsalgu@ucm.es}

\author[Tradacete]{Pedro Tradacete}
\address{Instituto de Ciencias Matem\'aticas (CSIC-UAM-UC3M-UCM)\\
Consejo Superior de Investigaciones Cient\'ificas\\
C/ Nicol\'as Cabrera, 13--15, Campus de Cantoblanco UAM\\
28049 Madrid, Spain.}
\email{pedro.tradacete@icmat.es}

\author[Trejo-Arroyo]{Nazaret Trejo-Arroyo}
\address{Departamento de An\'alisis Matem\'atico y Matem\'atica Aplicada\\
Universidad Complutense de Madrid\\
Plaza de las Ciencias 3\\
28040 Madrid, Spain.}
\email{ntrejo@ucm.es}

\date{\today}


\begin{document}
\begin{abstract} We study different versions of \emph{free objects} in the setting of quasi-Banach spaces and quasi-Banach lattices. Special attention is devoted to the free $p$-convex $p$-Banach lattice $\FpBLp[E]$ generated by a $p$-natural quasi-Banach space $E$, for which we provide a functional representation by means of operators into $L_p[0,1]$. This representation yields, among other consequences: (1) Operators from a Banach space $E$ to any $p$-convex $(0<p<1)$ quasi-Banach lattice $X$ can be extended to lattice homomorphisms $\FBL[E] \to X$ with control of the norm. (2) The space $\ell_p(\Gamma)$ $(0<p<1)$ is a projective $p$-Banach lattice precisely when $\Gamma$ is countable. 
(3) The free vector lattice generated by $E$ sits inside $\FpBLp[E]$ as a dense sublattice. 
\end{abstract}
	
    \maketitle    
	\section{Introduction}

\par The topics of $p$-Banach spaces and $p$-Banach lattices for $0 < p < 1$, although natural steps in extending the classical framework of functional analysis beyond locally convex settings, have not received as much attention as the well-established theories of Banach spaces and Banach lattices. This gap in our knowledge is usually attributed to the spectacular failure of Hahn-Banach's theorem in the absence of local convexity --see \cite[Sections 4 and 5]{f-space} for an account of the situation. However, there have been serious attempts to develop a satisfactory theory of $p$-Banach and quasi-Banach spaces, as it can be attested by the monographs \cite{aoki,f-space,rolewicz85} and the survey paper \cite{kalton-handbook}. 

\par It must be conceded that, apart from their intrinsic interest, quasi-Banach spaces are not an isolated subject with respect to the theory of Banach spaces. Indeed, it is not uncommon that a topic originating in the context of Banach spaces cannot be fully studied within it, or may benefit from a wider or more natural perspective when considered in a non-locally convex setting. Examples include, but are not limited to, harmonic analysis \cite{kalton-plurisub}, interpolation theory \cite{KM98} and exact sequences of Banach spaces \cite{kalton-3sp}. In  turn, quasi-Banach and $p$-Banach lattices provide a natural extension of Banach lattices and still permit the application of many classical techniques from their locally convex counterpart. In fact, several recurring themes present in Banach lattices have already been expanded in the absence of local convexity \cite{AAW21, kalton-conv, SP-T19}. Furthermore, since quasi-Banach lattices are often the most natural examples of quasi-Banach spaces, they are expected to play an important role in the subject, in a similar way as Banach lattices do in the theory of Banach spaces. 

\par On the other hand, the study of \emph{free objects} in the context of Banach lattices has developed a notable line of research. While the notion of free vector lattice over a vector space goes back to the 1960's \cite{baker,bleier}, the normed structure was first incorporated by B. de Pagter and A. Wickstead in \cite{depag-wick}, where they provided the \emph{free Banach lattice generated by a set}. This construction was the precursor of the \emph{free Banach lattice generated by a Banach space}, due to A. Avilés, J. Rodríguez and the second-named author \cite{ART}. Several variations of these theme have been subsequently addressed, with particular attention to $p$-convexity in Banach lattices \cite{fbl-convex,fbl}. 
Free Banach lattices have proven useful to explore the relationship between Banach space and Banach lattice properties, which is a recurring theme in the literature --see \cite{jmst, kalton-lattice}--, and recent contributions connect them to projectivity in Banach lattices \cite{AMR, AMR2}, size of disjoint families \cite{APR}, push-out constructions \cite{AT} and approximation properties \cite{Oikhberg}, among many others. There has also been further developments concerning free Banach lattices generated by a lattice \cite{AMRR, AR}, upper $p$-estimates \cite{GLTT} and free complex Banach lattices \cite{dHT}.

\par The present paper is concerned with the study of \emph{free quasi-Banach lattices} generated by a quasi-Banach space. Our motivation is to clarify the relationship between the developing theories of quasi-Banach spaces and quasi-Banach lattices, while at the same time providing a new tool that may contribute to the advancement of both. A second source of our motivation lies in earlier efforts in extending free objects from the Banach space to the quasi-Banach space setting, such as for the well-known class of \emph{Lipschitz free spaces} generated by a metric space. This constitutes a very active research field at the moment, and its $p$-Banach space analogue was first introduced in \cite{AK2009} and further studied in \cite{AACD,AACD2}.

\par We now outline the contents of the paper. The theory of quasi-Banach and $p$-Banach spaces suggests a number of natural routes to define free objects. Perhaps the most immediate is the \emph{free $p$-Banach lattice generated by a $p$-Banach space;} we devote Section \ref{sec:3} to prove this object indeed exists. However, it should be noted that the potential of free objects in functional analysis generally relies in having a concrete representation. The main obstruction in obtaining such representations in the non-locally convex setting is the absence of Hahn-Banach theorem, which played an essential role in the workable description of the free Banach lattice generated by a Banach space --given in \cite{ART}-- as a certain lattice of functions on the dual unit ball. To circumvent these difficulties, while also accounting for $p$-convexity conditions in quasi-Banach lattices --here $p>0$--, we restrict ourselves to the class of \emph{natural} quasi-Banach spaces, i.e. those admitting a linear embedding into a quasi-Banach lattice which is $p$-convex for some $p>0$. For these spaces, we develop a notion of ``dual'' space which allows us to prove not only existence, but a representation of the \emph{$p$-convex $p$-Banach lattice} as a certain space of bounded, positively homogeneous and continuous functions on a suitable topological space, and which can be seen as the right analogue in the non-locally convex setting. This is contained in Section \ref{sec:4}. 

\par The starting point of Section \ref{sec:5} is that $p$-Banach spaces are also $r$-Banach spaces for any $0<r<p$. Therefore, if a free $p$-convex $p$-Banach lattice exists for some $0<p\leq 1$, then there is a family of free quasi-Banach lattices parameterized by $r\in (0,p]$. The study of the relations between them leads us to a surprising conclusion: if $r<p$, then the free $r$-convex $r$-Banach lattice automatically extends operators to any quasi-Banach lattice which is $q$-convex for every $0<q<r$. This situation implies certain automatic convexity properties for the free object, in the spirit of \cite[Section 9.6]{fbl}, as well as relations between free $p$-Banach lattices and their natural counterparts in the category of quasi-Banach spaces. We also particularize the results of this section to the study of the \emph{free Banach lattice}, and show that they are optimal in the context of Banach spaces. Finally, in Section \ref{sec:6} we pursue the classical relation between free and projective objects by studying projectivity in the category of $p$-Banach lattices. Precisely, by considering chain conditions, we extend several results of \cite{AMRT22} by showing that, for $0<p<1$, the space $\ell_p(\Gamma)$ is a projective $p$-Banach lattice precisely when $\Gamma$ is countable. 
        
    \section{Preliminaries}
	\subsection{Quasi-norms and $p$-norms} Given $E$ a vector space, we say that a map $\|\cdot\|: E \to \R^{+}$ is a \emph{quasi-norm} if 
	\begin{enumerate}
		\item $\|e\|=0$ implies $e=0$,
		\item $\|\lambda e\| = |\lambda|\|e\|$ for every $e\in E$, $\lambda\in \R$,
		\item there exists $\Delta \geq 1$ so that $\|e+e'\|\leq \Delta(\|e\|+\|e'\|)$ for every $e,e'\in E$. 
	\end{enumerate}
	We also say that $\|\cdot\|$ is a \emph{$p$-norm}, for a fixed $p\in (0,1]$, if it satisfies (1), (2) and 
	\begin{enumerate}
		\item[(3')] $\|e+e'\|^p\leq \|e\|^p+\|e'\|^p$ for every $e,e'\in E$, 
	\end{enumerate}
	Quasi-norms (and $p$-norms) induce a linear topology on $E$ by means of the basis of neighborhoods at zero given by $B(0,\varepsilon) = \{ e: \|e\|<\varepsilon\}$ for $\varepsilon>0$. If such a topology is complete, then $(E, \|\cdot\|)$ is a \emph{quasi-Banach} space, or a \emph{$p$-Banach space}, respectively. Furthermore, if $E$ carries a lattice structure (i.e. a partial order with suprema $x\vee y$ and infima $x\wedge y$ for any pair of vectors $x,y\in E$) compatible with its quasi-norm (or $p$-norm) structure, in the sense that
	\begin{enumerate}
		\item[(4)] $x\leq y$ implies $x+z\leq y+z$, for every $z\in E$
		\item[(5)] $x\leq y$ implies $\lambda x \leq \lambda y$, provided $\lambda \geq 0$, 
		\item[(6)] $|x|\leq |y|$ implies $\|x\| \leq \|y\|$ (where $|x|:=x\vee (-x)$),
	\end{enumerate}
	then $E$ is a \emph{quasi-Banach lattice} (resp.,  \emph{$p$-Banach lattice}). A quasi-norm or a $p$-norm which satisfies (6) is called a lattice quasi-norm (resp., $p$-norm). Observe that every $p$-norm is a quasi-norm (for $\Delta = 2^{1/p-1}$), and the Aoki-Rolewicz theorem \cite[Thm 1.2]{f-space} asserts that every quasi-norm admits an equivalent $p$-norm for some $0<p\leq1$.	
	
	\subsection{Convexity conditions for quasi-Banach lattices }
	
	Given $p\in (0,+\infty)$, we say a quasi-Banach lattice $X$ is \emph{(lattice) $p$-convex} if there exists $C>0$ such that for every $x_1, \dots, x_n\in X$,
	\[ \left\| \left(\sum_{k=1}^n |x_k|^p\right)^\frac1p\right\| \leq C \left( \sum_{k=1}^n \|x_k\|^p \right)^\frac1p.  \]
    The above definition also works for the case $p=\infty$ if now one requires
    	\[ \left\| \bigvee_{k=1}^n |x_k|\right\| \leq C \left( \max_{1\leq k \leq n} \|x_k\|\right).  \]
    The infimum of the constants $C$ that satisfy the above inequality is referred to as the \emph{$p$-convexity constant} of $X$ and is denoted $M^{(p)}(X)$. Here, recall that positively homogeneous expressions involving finitely many variables, such as $\left(\sum_{k=1}^n |x_k|^p\right)^\frac1p$, are well defined elements in any quasi-Banach lattice by means of Yudin-Krivine functional calculus (cf. \cite[1.d]{lin-tza}). Furthermore, if $X$ is a $p$-convex quasi-Banach lattice, then, by a Hölder-type argument, $X$ is also $r$-convex for $0<r<p$, and $M^{(r)}(X) \leq M^{(p)}(X)$.
	
	\par Throughout this paper, the term \emph{``$p$-convex''} will always mean \emph{``lattice $p$-convex''}. We will employ the term \emph{``$p$-normable''} to refer to a quasi-Banach space whose quasi-norm is equivalent to a $p$-norm. It is clear that, given $p\in (0,1]$, every $p$-convex Banach lattice is $p$-normable (in fact it admits an equivalent lattice $p$-norm). The reciprocal is only true for $p=1$; indeed, for every $p\in (0,1)$, the space $L_{p,\infty}[0,1]$ is $p$-normable but it is $r$-convex only for $r<p$ \cite{Popa}. 
    
	\begin{defi} 
    A quasi-Banach lattice $X$ is \emph{L-convex} if there exists $\varepsilon\in(0,1)$ such that for any $u\in X_+$ with $\|u\|=1$ and $x_1, \dots, x_n \in [0,u]$ satisfying 
	\[ \frac{1}{n} \sum_{k=1}^n x_k \geq (1-\varepsilon)u,\]
	then $\max_{1\leq k \leq n} \|x_k\| \geq \varepsilon$.
    \end{defi} 
    
	The class of $L$-convex quasi-Banach lattices agrees with the class of quasi-Banach lattices which are $p$-convex for some $p>0$ \cite[\S2]{kalton-conv}. In particular, for every measure $\mu$, the spaces $L_p(\mu)$ are $L$-convex, and more generally, any Orlicz space whose Orlicz function satisfies the $\Delta_2$ condition. Perhaps the simplest example of a quasi-Banach lattice which is not $L$-convex is the space $L_p(\phi)$ of measurable real-valued functions on a set $\Omega$ such that 
    \[ \|f\|_p = \left(\int_{[0,+\infty)} \phi\{|f|>t^p\} dt\right)^\frac{1}{p} < +\infty,\]
    where $\phi$ is a \emph{pathological submeasure} (see e.g. \cite[Example 2.4]{kalton-conv} for details and concrete examples). 
    
	\par Given a quasi-Banach lattice $X$, and $p>0$, the \emph{$p$-convexification} of $X$ is the quasi-Banach lattice $X^{(p)}$ whose underlying vector lattice is $X$ with the operations 
	\[x \oplus y = (x^{1/p} + y^{1/p})^p, \qquad \lambda \odot x = \lambda^p x, \]
	and with the same order of $X$, and whose quasi-norm is $\|x\|_{(p)} = \|x\|^{1/p}$. It is easy to check that if $X$ is $r$-convex for some $r>0$, then $X^{(p)}$ is $pr$-convex. Furthermore, if $T: X \to Y$ is a lattice homomorphism between $p$-Banach lattices, one can check that $T: X^{(p)} \to Y^{(p)}$ is still a (linear) lattice homomorphism satisfying
        \[\|Tx\|_{(p)}  \leq \|T\|^{1/p} \|x\|_{(p)}. \]
    Similarly, one can check that $T:X\to Y$ is a lattice embedding precisely when $T:X^{(p)}\to Y^{(p)}$ is.

	\begin{defi} 
		Let $p\in (0,1]$. A quasi-Banach space $E$ is \emph{$p$-natural} if there exists a $p$-convex $p$-Banach lattice $X$ with $M^{(p)}(X)=1$ and a linear isometric embedding from $E$ to $X$. Furthermore, we say that a quasi-Banach space $E$ is \emph{natural} if it is $p$-natural for some $p\in (0,+\infty]$. 
	\end{defi}

    \par The previous definition is essentially due to Kalton \cite[p. 1118]{kalton-handbook} with a slight difference: we ask that the embedding from $E$ to $X$ is an \emph{isometry}. However, every quasi-Banach space $E$ which is isomorphic to a subspace of a $p$-convex $p$-Banach lattice admits an equivalent quasi-norm (in fact, a $p$-norm) which turns it into a $p$-natural quasi-Banach space. In particular, observe that if $p \in (0, 1]$, then every $p$-natural quasi-Banach space \emph{is} automatically a $p$-Banach space (here ``$1$-Banach space'' just means ``Banach space''). For $p=1$ we can say more: every Banach space $E$ embeds isometrically into $C(B_{E^*})$, which is an $\infty$-convex Banach lattice.

    \par Natural spaces can be regarded as a relatively well-behaved class of quasi-Banach spaces; and in fact, a large number of quasi-Banach spaces appearing in analysis are natural. Nevertheless, one can find several non-natural quasi-Banach spaces: $L_p(\mathbb T)/H_p(\mathbb T)$ \cite[Section 10]{kalton-linear}, the Schatten classes for $0<p<1$ \cite[Corollary 4.9]{kalton-plurisub}, and $L_p(\phi)$ for $\phi$ any pathological submeasure \cite[Example 2.4]{kalton-conv} constitute some well-known examples. In addition, there is a remarkable example of a non-natural quasi-Banach space without any basic sequence \cite{kalton-basis}.

	 \subsection{Free Banach lattices}\label{ssec:2.3}
     Given a Banach space $E$ and $p\in [1,+\infty]$, the \emph{free $p$-convex Banach lattice} generated by $E$ is a $p$-convex Banach lattice $\FBLp[E]$ with $M^{(p)}(\FBLp[E])=1$, together with a linear isometric embedding $\delta: E \to \FBLp[E]$, such that every operator $T$ from $E$ into a $p$-convex Banach lattice $X$ admits a unique lattice homomorphism $\widehat{T}: \FBLp[E] \to X$ such that the following diagram commutes  
        \[ \begin{tikzcd}
            \FBLp[E]  \arrow[rd, "\widehat{T}"] & \\
            E \arrow[r, "T"] \arrow[u, "\delta"] & X
        \end{tikzcd}
     \]
     and with $\|\widehat{T}\|\leq M^{(p)}(X)\|T\|$. Since every Banach lattice is 1-convex with $M^{(1)}(X)=1$, the case $p=1$ provides the \emph{free Banach lattice}. 
     
     \par The construction of the free Banach lattice $\FBL[E]$ generated by a Banach space $E$ is developed in \cite{ART}, and later extended to the $p$-convex setting in \cite{fbl-convex}. Let us recall the key ideas of the explicit functional representation of $\FBLp[E]$. Denote by $H[E]$ the linear subspace of $\mathbb{R}^{E^*}$ consisting of all positively homogeneous functions $f:E^* \to \mathbb R$. Given $p\in [1,+\infty)$, every operator $T: E \to \ell_p^n$ is of the form $T(x) = \sum_{k=1}^n x^*_k(x)$, and hence admits a natural extension $\widehat{T}: H[E] \to \ell_p^n$ given by $T(f) = (f(x_k^*))_{k=1}^n$. This motivates the following definition for $f \in H[E]$: 
     \begin{align*}
         \notag
         \norm{f}_{\FBLp[E]}= & \sup \left\{ \left( \sum_{i=1}^n |f(x_k^*)|^p \right)^{1/p} \colon n\in\mathbb N, (x_i^*)_{i=1}^n \subset E^*, \sup_{x \in B_{E}} \sum_{i=1}^n |x_i^*(x)|^p \leq 1 \right\}  \\[1mm]
         =& \sup \left\{ \|\widehat{Tf}\|_{\ell_p^n}: n\in\mathbb N, T\in B_{\mathcal L(E, \ell_p^n)}\right\}. 
     \end{align*}
     Observe that, if for a given $x \in E$, we let $\delta_x \in H[E]$ be defined by $\delta_x(x^*)=x^*(x)$ for any $x^* \in E^*$, then $\norm{\delta_x}_{\FBLp[E]}=\norm{x}$ for every $x \in E$. Now, consider \[ H_p[E]=\lbrace f \in H[E]: \norm{f}_{\FBLp[E]} < +\infty \rbrace\] equipped with the pointwise linear and lattice operations, which becomes a $p$-convex Banach lattice under the norm $\norm{\cdot}_{\FBLp[E]}$. The space $\FBLp[E]$ can be identified with the closed sublattice generated by $\{\delta_x: x\in E\}$ in $H_p[E]$, together with the map $\delta:E \to \FBLp[E]$ given by $\delta(x)=\delta_x$. Indeed, the content of \cite[Theorem 2.1]{fbl} is precisely that this sublattice of $H_p[E]$ has the universal property of extension of operators on $E$ to lattice homomorphisms via $\delta$ (see also \cite[Theorem 2.5]{ART} and \cite[Theorem 6.1]{fbl-convex}).
     For the $\infty$-convex case, it is shown in \cite[Proposition 2.2]{fbl} that $\FBL^{(\infty)}[E]$ coincides with the lattice $C_{ph}(B_{E^*})$ of positively homogeneous weak$^*$ continuous functions on $B_{E^*}$.
     
\section{The free $p$-Banach lattice}\label{sec:3}
    The existence of a dual space is fundamental for the explicit construction of $\FBLp[E]$ in the Banach space setting, as it has been made clear in the previous section. Our goal now is to obtain the \emph{free object} in the context of $p$-Banach spaces, where  the existence of sufficient functionals is no longer guaranteed. Therefore, it is clear that we need to make use of different techniques. 
	\begin{defi}
		Consider $E$ a $p$-Banach space. The \emph{free $p$-Banach lattice over $E$} is a $p$-Banach lattice $\FpBL[E]$, together with an isometric embedding $\delta: E \hookrightarrow \FpBL[E]$, such that for every operator $T: E \to X$, where $X$ is a $p$-Banach lattice, there exists a unique lattice homomorphism $\widehat{T}: \FpBL[E] \to X$ \emph{extending $T$} --i.e., $\widehat{T} \delta = T$-- and such that $\|\widehat{T}\| = \|T\|$.
	\end{defi}
	\par Observe that the universal property of $\FpBL[E]$ precisely means the existence of a \emph{natural correspondence} 
    \[ \begin{array}{rcl} \mathcal L(E, X) & =\joinrel= & \Hom(\FpBL[E], X), \\ T & \longmapsto & \widehat{T} \end{array} \]
    where $\Hom(\FpBL[E], X)$ denotes the set of all lattice homomorphisms between $\FpBL[E]$ and any fixed $p$-Banach lattice $X$. 
    This property shows that, if $\FpBL[E]$ exists, then it is unique up to a lattice isomorphism. We now show by an abstract argument that $\FpBL[E]$ does exists for every $p$-Banach space $E$. The proof is essentially a combination of standard arguments from universal algebra (see also \cite{de-jeu}) and \cite{turpin}. 
	
	\begin{thm} \label{thm:FpBL}
		For every $p$-Banach space $E$, $\FpBL[E]$ exists. 
	\end{thm}
	\begin{proof} Let $\FVL[E]$ be the \emph{free vector lattice} associated to $E$. Such space comes equipped with an injective linear map $\delta: E \to \FVL[E]$ and the following \emph{universal property}: given a vector lattice $X$ and a linear map $T:E \to X$, there exists a unique linear homomorphism $\hat{T}:\FVL[E] \to X$ such that $\hat T\delta = T$.
    Consider the following $p$-seminorm on $\FVL[E]$:
		$$\nnorm{f}=\sup \|\hat{T}f\|,$$
		where the supremum runs over all operators $T: E \to X$ such that $\|T\|=1$ and $X$ is any $p$-Banach lattice. Actually, it is enough to consider operators into $p$-Banach lattices of density character not bigger than that of $E$, to avoid set theoretical inconsistencies. 
		It is easy to check that $\nnorm{\cdot}$ is a lattice $p$-seminorm. To make it a lattice $p$-norm we consider the ideal $\mathcal{N}=\lbrace f \in \FVL[E]: \nnorm{f}=0 \rbrace, $
		and define $\FpBL[E]$ as the completion of the quotient of $(\FVL[E], \nnorm{\cdot})$ by $\mathcal{N}$. 

		The space $\FpBL[E]$ is a $p$-Banach lattice endowed with the obvious $p$-norm, and the linear map $\delta: E \hookrightarrow \FpBL[E]$ becomes an isometric embedding, as we now show. The inequality $\nnorm{\overline{\delta_e}}\leq \norm{e}$ is clear by definition. For the converse, observe that, thanks to \cite[Proposition 2]{turpin}, there exist a $p$-Banach lattice $X$ and an isometric embedding $\iota: E \to X$. Therefore,  
		\[\norm{e}=\norm{\iota(e)}=\norm{\hat{\iota}(\delta_e)}\leq \nnorm{\overline{\delta_e}}=\nnorm{\delta_e}.\qedhere\]  
	\end{proof}

    \begin{remark} The reader may be prompted to think that the embedding $\iota$ mentioned in the last lines of Theorem \ref{thm:FpBL} is a mere abstract tool. This is not the case; in fact, a suitable modification of the proof of \cite[Proposition 2]{turpin} provides another proof of the existence of $\FpBL[E]$. We sketch the argument next: Given $E$ any $p$-Banach space, let $E^{\sharp}$ be the algebraic dual space of $E$, and consider the canonical inclusion $\delta: E \to \mathbb{R}^{E^\sharp}$ taking any $e\in E$ to the linear map $\delta_e$ given by the formula $\delta_e(e^\sharp) = e^\sharp(e)$. Then, the sublattice $E_0$ generated by $\{\delta_e: e\in E\}$ is lattice and linear isomorphic to $\FVL[E]$; in particular, every operator $T$ from $E$ to  any $p$-Banach lattice $X$ admits a unique lattice homomorphism $\widehat{T}: E_0 \to X$ such that $T = \widehat{T}\delta$.
    \par We now consider the following lattice $p$-seminorm on $E_0$:
	\[ |f|_0 = \inf \left\{ \left(\sum_{k=1}^n \|e_k\|^p \right)^\frac1p : |f| \leq \sum_{k=1}^n |\delta_{e_k}| \right\} \]
    Then, the inclusion $\delta: E \to E_0$ becomes an isometric embedding, just as in the proof of \cite[Proposition 2]{turpin}. The crux of the argument is that $|\cdot|_0$ also makes continuous all extensions $\widehat{T}$, since if $f\in E_0$ satisfies $|f| \leq \sum_{k=1}^n |\delta_{e_k}|$ for certain $e_1, \dots, e_n\in E$, then 
    \begin{align*} \|\widehat{T}f\| & \leq \left\|\sum_{k=1}^n\widehat{T}(\delta_{e_k})\right\| = \left\|\sum_{k=1}^n T(e_k)\right\| \leq \left(\sum_{k=1}^n \|T(e_k)\|^p\right)^{\frac1p} \leq \\[2mm]
    & \leq \|T\|\cdot \left(\sum_{k=1} \|e_k\|^p \right)^\frac1p 
    \end{align*}
    Therefore, if $\mathcal N=\{f\in E_0: |f|_0 = 0\}$, then the universal property of $\FpBL[E]$ asserts that the completion of $E_0/\mathcal N$ is isometrically and lattice isomorphic to $\FpBL[E]$. 
    \end{remark}
	
	\section{The free $p$-convex $p$-Banach lattice}\label{sec:4}
	\subsection{Representation of $\FpBL^{(p)}[E]$.}\label{section:representation FpBLp}
    Once the existence of the free object has been proven in an abstract setting, we will specialize to the category of $p$-natural quasi-Banach spaces. In this context, we give a functional representation of the \emph{free $p$-convex $p$-Banach lattice}. The key lies on the fact that, for a $p$-natural quasi-Banach space $E$, there are plenty of operators $E \to L_p[0,1]$ --see Proposition \ref{prop:p-natural} below--. 
    
    \par Indeed, fix $p\in (0,1]$ and let $E$ be a $p$-natural quasi-Banach space. Let us consider the ``\emph{$\vartriangle$-dual space}'', $E^{\vartriangle}=\mathcal{L}(E, L_p[0,1])$, equipped with $\|e^{\vartriangle}\|=\sup_{\|e\|\leq1}\|e^{\vartriangle}(e)\|_p$ for $e^{\vartriangle}\in E^{\vartriangle}$. This is clearly a $p$-Banach space and every $e \in E$ defines an operator $\delta_e: E^\vartriangle \to L_p[0,1]$ by the formula $\delta_e(e^{\vartriangle})=e^{\vartriangle}(e)$. Therefore, we will call $E^{\vartriangle \vartriangle}=\mathcal{L}(E^{\vartriangle}, L_p[0,1])$ the \emph{$\vartriangle$-bidual} of $E$. On the other hand, the unit ball $B_{E^{\vartriangle}}$ of $E^{\vartriangle}$ can be endowed with the \emph{``weak-$\vartriangle$'' topology} $\sigma(E^{\vartriangle},E)$ which is the coarsest topology which makes continuous the maps $\lbrace \delta_e: e \in E \rbrace$. 
    
    \par Of course, these notions can be very well vacuous for a general quasi-Banach space; after all, there exist quasi-Banach spaces which do not admit nonzero operators to any $L_p(\mu)$-space \cite[Theorem 10.3]{kalton-linear}. However, if $E$ is $p$-natural for some $p\in (0,1)$, we now show $E^\vartriangle$ behaves in many aspects like an authentic dual space of $E$.
    
	\begin{prop}\label{prop:p-natural}
		For a quasi-Banach space $E$,  the following assertions are equivalent:
		\begin{itemize}
			\item[\textit{i)}] $E$ is $p$-natural. 
			\item[\textit{ii)}] There exists a set $J$ with $|J|=\dens(E)$ and a linear isometric embedding $\iota:E \hookrightarrow \ell_{\infty}(J, L_p[0,1])$.
            \item[\textit{iii)}] The map $\delta: E \to E^{\vartriangle \vartriangle}$ is a linear isometric embedding.
            \item[\textit{iv)}] There exists a linear isometric embedding $E \hookrightarrow C_{ph}^b(B_{E^{\vartriangle}},L_p[0,1])$, where $C_{ph}^b(B_{E^\vartriangle}, L_p[0,1])$ is the space of bounded positively homogeneous continuous functions $f: (B_{E^\vartriangle}, \sigma(E^\vartriangle, E)) \to L_p[0,1]$.
		\end{itemize}
	\end{prop}
	\begin{proof} $(i) \Rightarrow (ii)$.-- Assume $E$ isometrically embeds in a $p$-convex quasi-Banach lattice $X$, and consider $X^{(1/p)}$ the $(1/p)$-convexification of $X$, which is 1-convex by hypothesis and therefore normable. 
    Furthermore, if $\kappa=\dens(E)$, we can assume that $\dens(X)=\dens(X^{(1/p)})=\kappa$ just replacing $X$ by the sublattice generated by $E$ inside $X$. 
    Let $I$ be a set with $|I|=\kappa$ and choose $\lbrace x_i\rbrace_{i \in I}$ a dense set in the unit ball of $X^{(1/p)}$. For each $i\in I$, choose $x_i^* \in (X_+^{(1/p)})^*$ with $\norm{x_i^*}=1$ such that $x_i^*|x_i|=\norm{x_i}$. Then, setting $N_i=\lbrace x \in B_{X^{(1/p)}}:x_i^*|x|=0\rbrace$, by Kakutani's representation theorem the completion of the quotient space $(X^{(1/p)}/N_i, x_i^*(|\cdot|))$ is isomorphic to some $L_1$-space, say $L_1(\mu_{x_i^*})$. The canonical quotient map induces a lattice homomorphism  $\pi_i: X^{(1/p)} \to L_1(\mu_{x_i^*})$ with dense range. Consider now the map 
		\[
		\Pi:X^{(1/p)} \longrightarrow \left(\oplus_{i\in I}(L_1(\mu_{x^*_i})\right)_\infty, \qquad x \longmapsto \Pi(x):=(\pi_i(x))_{i \in I}
		\]
		which is a lattice isometric embedding, for if $x\in X^{(1/p)}$ with $\norm{x}\leq 1$ and $\varepsilon>0$ there is $i\in I$ such that $\norm{x-x_i}\leq \varepsilon$, and so,
		\begin{align*}
			|x^*_i(x)| & >\big||x^*_i(x-x_i)|-|x^*_i(x_i)|\big| \geq \norm{x_i}-\varepsilon \geq \norm{x}-\norm{x-x_i}-\varepsilon\geq \\ & \geq \norm{x}-2\varepsilon.
		\end{align*}
		Hence, $\norm{\pi_i(x)}_{L_1(\mu_{x_i^*})}\geq \norm{x}-2\varepsilon$, so $\|\Pi(x)\|_\infty = \|x\|$. 

        \begin{claim} $L_1(\mu_{x_i^*})$ is isometrically and lattice isomorphic to a sublattice of \[ L_1(M_I):= \ell_1(I)\oplus_1 L_1([0,1]^{I}).\]   
        \end{claim}
        \begin{proof}[Proof of the Claim.]\renewcommand{\qedsymbol}{$\blacksquare$} 
        Since $\pi_i$ has dense range,  $\dens [L_1(\mu_{x_i^*})] \leq \kappa$, and so, an appeal to Maharam's theorem \cite[\S14, Theorems 9 and 10]{lacey} yields that $L_1(\mu_{x_i^*})$ is isometrically and lattice isomorphic to a sublattice of $\ell_1(I)\oplus_1 \ell_1(L_1([0,1]^{I}))$.
        Now, it is well-known that $\ell_1(L_1[0,1]^I)$ is isometrically and lattice isomorphic to a sublattice of $L_1([0,1]^I)$. Indeed, fix $i_0\in I$ and consider the following sequence of disjuoint subspaces of $[0,1]^I$: 
        \[ \Omega_n = \{x\in [0,1]^I : 2^{-n} \leq x(i_0) < 2^{-n+1} \}. \]
        Observe that $\mu(\Omega_n) = 2^{-n}$ and that every $\Omega_n$ is homeomorphic to $[0,1]^I$.  Hence, if $\psi_n: \Omega_n \to [0,1]^I$ are the corresponding homeomorphisms, then the map 
        \[\ell_1(L_1([0,1]^I)) \longrightarrow L_1([0,1]^I), \quad (g_n)_{n=1}^\infty \longmapsto g(x) = \sum_{n=1}^\infty 2^n\cdot g_n(\psi_n(x))\cdot I_{\Omega_n}(x), \]
        is a lattice isometry. 
        \end{proof}

        All in all, $\Pi$ induces a lattice isometric embedding 
		\begin{equation} \label{eq:natural-1} X^{(1/p)}\overset{\Pi}{\longrightarrow} \left(\oplus_{i\in I}(L_1(\mu_{x^*_i})\right)_\infty \hookrightarrow \ell_\infty\big(I, L_1(M_I)\big).\end{equation}
		We now need to show that  $L_1(M_I)$ is isometrically and lattice isomorphic to a sublattice of $\ell_\infty(J, L_1[0,1])$, where $J$ is some set of indices with $|J| = \kappa$. Let us denote $\Fin(I)$ for the collection of finite subsets of $I$. Fix $G \in \Fin(I)$ and consider the map 
		\[ \iota_G: [0,1]^G \longrightarrow [0,1]^I, \qquad \iota_G(x)(i) = \begin{cases} x(i) & \text{ if } i\in G, \\ 0 & \text{ if } i \in I \setminus G. \end{cases} \]
		Then, given a pair $(F,G)$ of finite subsets of $I$, the continuous map
		\[\rho_{(F,G)}: L_1(M_I) \longrightarrow \ell_1(F) \oplus_1  L_1([0,1]^G), \qquad \rho_{(F,G)}(f,g) = (\restr{f}{F}, g\circ \iota_G)\]
		is a contractive lattice homomorphism.

        \begin{claim}
            For every $(f,g)\in L_1(M_I)$, the following equality holds:
            \[ \|(f,g)\|_{L_1(M_I)} = \sup_{(F,G)\in \Fin(I)^2} \|\rho_{(F,G)}(f,g)\|_{\ell_1(F) \oplus_1  L_1([0,1]^G)} \]
        \end{claim}
        \begin{proof}[Proof of the Claim.]\renewcommand{\qedsymbol}{$\blacksquare$} Given $(f,g)\in L_1(M_I)$ and $\varepsilon>0$, there is $F \in \Fin(I)$ such that $\|f\|_{\ell_1(I)} \leq \|\restr{f}{F}\|_{\ell_1(F)} + \varepsilon$. 
        To deal with $g\in L_1([0,1]^I)$, we observe that the set of continuous functions from $[0,1]^I$ into $[0,1]$ which depend on finitely many coordinates is dense in $L_1([0,1]^I)$ (just apply the Stone-Weierstrass theorem). 
        Hence there exists $g_\varepsilon \in L_1([0,1]^I)$ depending only on the coordinates of a finite subset $G \subset I$ and such that $\|g_\varepsilon - g\|_{L_1[0,1]^I}<\varepsilon$. 
        Then, by virtue of Fubini's theorem  we have that $\|g_\varepsilon\|_{L_1([0,1]^I)} = \|g_\varepsilon \circ \iota_G\|_{L_1([0,1]^G)}$, and so $\|g\|_{L_1([0,1])^I} \leq \|g\circ \iota_G\|_{L_1([0,1]^G)} + 2\varepsilon$. Combining this, we get that 
		\[\|f\|_{\ell_1(I)} + \|g\|_{L_1([0,1]^I} \leq \|\restr fF\|_{\ell_1(F)} + \|g\circ \iota_G\|_{L_1([0,1]^G)} + 3\varepsilon.\]
		This is enough to conclude. \end{proof}
                As a consequence, the map 
		\begin{gather*} R: L_1(M_I) \longrightarrow  \oplus_{(F,G)} \left( \ell_1(F) \oplus_1  L_1([0,1]^G)\right)_\infty,  \\  R(f,g) = (\rho_{(F,G)}(f,g))_{(F,G)}, 
        \end{gather*}
		is a lattice isometric embedding.       
        Finally, observe that for every $(F,G)\in \Fin(I)^2$, $\ell_1(F) \oplus_1 L_1([0,1]^G)$ is a separable $L_1$-space, hence it is isometrically lattice isomorphic to a sublattice of $L_1[0,1]$. This induces an isometric lattice embedding 
		\begin{equation} \label{eq:natural-2} \begin{aligned} \ell_\infty\left(I, L_1(M_I) \right) \overset{\oplus_{i\in I} R}{\xrightarrow{\hspace{1cm}}} & \ \ell_\infty\left(I, \oplus_{(F,G)} \left( \ell_1(F) \oplus_1 L_1([0,1]^G)\right) \right)_\infty  \\[1mm]
				\hookrightarrow &\ \ell_\infty\left(I, \ \ell_\infty\!\left(\Fin(I)^2,L_1[0,1]\right)_\infty\right)_\infty \\[1mm]
				= & \ \ell_\infty(J, L_1[0,1]),	
		\end{aligned} \end{equation}
		where $J=I \times \Fin(I)^2$ has the cardinality of $I$. 
		\par To finish, we combine Equations \ref{eq:natural-1} and \ref{eq:natural-2} to produce a lattice isometric embedding $T: X^{(1/p)}  \to \ell_\infty(J, L_1[0,1])$. Since $E$ embeds isometrically into $X$ by hypothesis, the composition of such an embedding $E \hookrightarrow X$ with the $p$-convexification of $T$ provides the desired isometric embedding from $E$ into $\ell_\infty(J, L_p[0,1])$.
        \par $(ii) \Rightarrow (iii)$.-- If there exists an isometric embedding $\iota: E \to \ell_\infty(J, L_p[0,1])$, then $E^\vartriangle \neq 0$; in fact, for every $j\in J$, we have a norm-one operator $T_j: E \to L_p[0,1]$ given by $T_j(e) = \iota(e)(j)$. This implies that $\delta: E \longrightarrow E^{\vartriangle\vartriangle}$ is an injective norm-one operator. Furthermore, 
        \[ \|\delta_e\| = \sup_{e^{\vartriangle}\in B_{E^\vartriangle}} \|e^{\vartriangle}(e)\| \geq \sup_{j\in J} \|T_j(e)\| = \|\iota(e)\| = \|e\|,\]
        that is to say, $\delta$ is an isometric embedding. 
        
        \par $(iii) \Rightarrow (iv)$.-- Just observe that the image of $\delta: E \longrightarrow E^{\vartriangle\vartriangle}$ actually lies in $C_{ph}^b(B_{E^{\vartriangle}},L_p[0,1])$. 

        \par $(iv) \Rightarrow (i)$.-- It is straightforward to check that $C_{ph}^b(B_{E^{\vartriangle}},L_p[0,1])$ with the pointwise order and lattice operations equipped with the $p$-norm $\|f\|=\sup_{e^{\vartriangle}\in B_{E^{\vartriangle}}}\|f(e^{\vartriangle})\|_p$ is a $p$-Banach lattice. Since $L_p[0,1]$ is $p$-convex, for every $e^{\vartriangle}\in B_{E^\vartriangle}$ we have
		\[
		\norm{\left(\sum_{k=1}^n |f_k(e^{\vartriangle})|^p\right)^{1/p}}_{L_p[0,1]} = \left( \sum_{k=1}^n \norm{f_k(e^{\vartriangle})}^p_{L_p[0,1]}\right)^{1/p}\leq \left( \sum_{k=1}^n \norm{f_k}_{\infty}^p\right)^{1/p}
		\]
		and therefore:
		\begin{align*}
			\norm{\left( \sum_{k=1}^n |f_k|^p \right)^{1/p}}_{\infty} & =
			\sup_{e^{\vartriangle}\in B_{E^\vartriangle}}\norm{\left[ \left( \sum_{k=1}^n |f_k|^p\right)^{1/p}\right](e^{\vartriangle})}_{L_p[0,1]}
			\\[2mm]
			& =\sup_{e^{\vartriangle}\in B_{E^\vartriangle}} \norm{\left( \sum_{k=1}^n |f_k(e^{\vartriangle})|^p \right)^{1/p}}_{L_p[0,1]} \leq \left( \sum_{k=1}^n \norm{f_k}_{\infty}^p\right)^{1/p}. 
		\end{align*}
        Hence, $C_{ph}^b(B_{E^{\vartriangle}},L_p[0,1])$ is a $p$-convex quasi-Banach lattice where $E$ embeds isometrically.
	\end{proof}
	
	\begin{defi} \label{defi:FpBLp} Let $E$ be a $p$-natural quasi-Banach space. The \emph{free $p$-convex $p$-Banach lattice} generated by $E$ is a $p$-convex $p$-Banach lattice $\FpBLp[E]$ with $M^{(p)}(\FpBLp[E]) = 1$, together with an isometric embedding $\delta: E \to \FpBLp[E]$, such that for every operator $T:E\to X$, where $X$ is a $p$-convex $p$-Banach lattice, there is a unique lattice homomorphism $\widehat{T}: \FpBLp[E] \to X$ such that $\widehat{T}\delta = T$ and $\|\widehat{T}\| \leq M^{(p)}(X) \|T\|$. \end{defi}
	
	\par Taking $p=1$ in the definition of $\FpBLp[E]$ recovers the  \emph{free Banach lattice generated by a Banach space $E$}, which is denoted $\FBL[E]$. Indeed, observe that $1$-natural quasi-Banach spaces are exactly Banach spaces, and that every Banach lattice $X$ is $1$-convex with $M^{(1)}(X) = 1$. 
	
	\par As usual, the universal property of $\FpBLp[E]$ implies uniqueness, provided it exists. Also note that if $\FpBLp[E]$ exists, then as $\delta$ defines a linear isometric embedding of $E$ into $\FpBLp[E]$, necessarily $E$ must be $p$-natural. We now show that this is the only requirement for the existence of $\FpBLp[E]$. 
	
	\begin{thm} 
		For every $p$-natural quasi-Banach space $E$, $\FpBLp[E]$ exists. 
	\end{thm}
	\begin{proof}
		
        \par Let $E$ be a $p$-natural quasi-Banach space and let us denote $\delta: E \to C_{ph}^b(B_{E^\vartriangle}, L_p[0,1])$ for the canonical embedding introduced in Proposition \ref{prop:p-natural}. Observe that, for every $T \in B_{E^\vartriangle}$, there is an extension
		\[ \widehat{T}:C_{ph}^b(B_{E^\vartriangle}, L_p[0,1]) \longrightarrow L_p[0,1], \qquad f\longmapsto \widehat{T}(f):=f(T)\]
		with $\|\widehat{T}\|=\|T\|$. 
        Hence, we will declare $\FpBLp[E]$ to be the closed sublattice of $C_{ph}^b(B_{E^\vartriangle}, L_p[0,1])$ generated by $\{\delta_e: e\in E\}$ and check that it satisfies Definition \ref{defi:FpBLp}. First of all, $\FpBLp[E]$ is $p$-convex with $p$-convexity constant equal to 1, since it is a closed sublattice of $C_{ph}^b(B_{E^\vartriangle}, L_p[0,1])$.  

		\par All that is left is to show that $\FpBLp[E]$ enjoys the \emph{universal property}. Let $X$ be a $p$-convex $p$-Banach lattice with $M^{(p)}(X)=1$ and let $T:E \to X$ be an operator with $\norm{T}=1$. Without loss of generality, we can assume that $\dens(X)=\dens(E)=\kappa$ (replace $X$ with the closed sublattice of $X$ generated by $T(E)$). Consider, by appealing to Proposition \ref{prop:p-natural}, a lattice isometric embbeding $\rho: X \hookrightarrow \ell_{\infty}(I, L_p[0,1])$, where $|I|=\kappa$, and write $\pi_i: \ell_{\infty}(I, L_p[0,1]) \to L_p[0,1]$ for the canonical projection defined by a given $i\in I$. Then, the maps  
		\[ T_{i}: E \to L_p[0,1], \qquad T_{i} = \pi_{i}\rho T,\] 
		admit extensions $\widehat{T_{i}}: \FpBLp[E] \to L_p[0,1]$ with $\|\widehat{T_{i}}\|= \norm{T_{i}} \leq \|T\|$.
        
		Note that, if $f\in \FpBLp[E]$ is of the form $f=\bigvee_{j=1}^n \delta_{e_j}-\bigvee_{k=1}^m \delta_{e'_k}$, where $e_1, \dots, e_n, e'_1, \dots, e'_m\in E$, then $(\widehat{T_{i}}(f))_{i\in I} \in \ell_\infty(I, L_p[0,1])$ actually lies in $\rho(X)$. Indeed, this follows from the identity
		\[ \widehat{T_{i}}(f) = \pi_i\rho\left(\bigvee_{j=1}^n T(e_j)-\bigvee_{k=1}^m T(e'_k)\right),\]
		for $i\in I$. This allows us to define $\widehat T(f)=\rho^{-1}((\widehat{T_{i}}(f))_{i\in I})\in X$.

        Now, it follows that
        \begin{equation*}
			\|\hat{T}f\|_X=\|\rho(\hat{T}f)\|_{\ell_\infty(I, L_p[0,1])}= \sup_{i\in I} \| \widehat{T}_{i}(f)\|_{L_p[0,1]}  \leq  \norm{T} \norm{f}_{\FpBLp[E]}. 
		\end{equation*} 
        Since every element of $\FpBLp[E]$ can be approximated with respect to $\|\cdot\|_{\FpBLp[E]}$ by elements of the form $f=\bigvee_{j=1}^n \delta_{e_j}-\bigvee_{k=1}^m \delta_{e'_k}$, the lattice homomorphism $\widehat{T}:\FpBLp[E]\rightarrow X$ is well-defined, unique and satisfies $\|\widehat{T}\|=\|T\|$ as desired.		
	\end{proof}

	\begin{remark}  
		If $E$ is an infinite-dimensional $p$-natural quasi-Banach space, for $0<p\leq 1$, then there is no compact Hausdorff vector topology on $B_{E^\vartriangle}$, and therefore the above representation does not allow to view the elements of $\FpBLp[E]$ as continuous functions on a \emph{compact} space. Indeed, the main result in \cite{kalton-note} asserts that, there is no Hausdorff vector topology on $L_p[0,1]$ that makes its unit ball relatively compact. Now, consider $\tau$ and $\tau'$ any two Hausdorff vector topologies in $B_{E^\vartriangle}$ and $L_p[0,1]$, respectively, such that all evaluation maps $\delta_e: (B_{E^\vartriangle}, \tau) \to (L_p[0,1], \tau')$ are continuous. Since  $L_p[0,1]$ is a transitive space --see \cite[Theorem 7.4]{f-space} or next subsection--, we infer that $\delta_e$ maps $B_{E^\vartriangle}$ onto the unit ball of $L_p[0,1]$ for any norm-one $e\in E$. Therefore, compactness of $(B_{E^\vartriangle}, \tau)$ would imply compactness of $(B_{L_p[0,1]}, \tau')$, which is a contradiction. This is in clear contrast with the representation of $\FpBL[E]$ for $p\geq1$ when $E$ is a Banach space given in \cite{fbl-convex}.
	\end{remark}

	\subsection{The structure of $p$-natural quasi-Banach spaces}

    \par We dedicate a moment to collect several properties of $p$-natural quasi-Banach which will be later needed for the study of free $p$-convex $p$-Banach lattices.

    \begin{defi} A quasi-Banach space $E$ is \emph{strictly transitive} if for every collection of linearly independent vectors $x_1, \ldots, x_n\in E$ and every collection $y_1, \ldots, y_n\in E$, there is $T\in \mathcal L(E,E)$ such that $T(x_i)=y_i$ for all $i\in \{1, \ldots, n\}$.
    \end{defi}

    \par The space $L_p[0,1]$ is strictly transitive \cite[\S2]{kalton-quotient} --cf. \cite[Theorem 7.5]{f-space}--. We will make use of this fact several times in the arguments below. Also, every quasi-Banach space with separating dual is strictly transitive (see Remark \ref{rem:transitive} for a proof).

    \par The main result of this section is the following: 

     \begin{thm}\label{th:surjective}
     Let $E$ be a $p$-natural quasi-Banach space. If $\{e_1,\ldots,e_n\}$ is a finite linearly independent subset of $E$, then the map 
     \begin{equation}
            \label{eq:epi_map}
        \begin{array}{rccl}        
            (e_i)_{i=1}^n \colon & E^{\vartriangle} & \longrightarrow & L_p[0,1] \times \overset{(n)}{\ldots} \times L_p[0,1] \\
             & e^{\vartriangle} & \longmapsto & (e^{\vartriangle}(e_i))_{i=1}^n
        \end{array}
        \end{equation}
     is surjective. 
     \end{thm}
     The proof proceeds by induction on $n$. The case $n=1$ is an immediate consequence of the fact that $L_p[0,1]$ is (strictly) transitive. Hence we treat the case $n=2$ as our initial case, which is in fact the crucial part of the proof. To prepare for the proof we need the following observation: 
    \begin{lemma}\label{lem:no0}
        Let $E$ be a $p$-natural quasi-Banach space. Given linearly independent $e_1, e_2 \in E$, there exists $e^{\vartriangle} \in E^{\vartriangle}$ such that $e^{\vartriangle}(e_1) \neq 0$ and $e^{\vartriangle}(e_2) \neq 0$.
    \end{lemma}
    \begin{proof}
        Since $E$ is $p$-natural, we can find $e^{\vartriangle} \in E^{\vartriangle}$ be such that $e^{\vartriangle}(e_1)\neq 0$. If $e^\vartriangle(e_2)\neq 0$ as well, then we are done. Otherwise, let us take $\widetilde e^{\vartriangle} \in E^{\vartriangle}$ such that $\widetilde e^{\vartriangle}(e_2) \neq 0$. Observe that for every $\lambda\in\R\backslash\{0\}$ we have
        \begin{itemize}
            \item $(e^{\vartriangle}+\lambda \widetilde e^{\vartriangle})(e_2)=\lambda \widetilde e^{\vartriangle}(e_2) \neq 0$, 
            \item $(e^{\vartriangle}+\lambda \widetilde e^{\vartriangle})(e_1)=e^{\vartriangle}(e_1)+\lambda \widetilde e^{\vartriangle}(e_1)$, 
        \end{itemize}
        and we can always choose $\lambda$ appropriately so that $e^{\vartriangle}(e_1)+\lambda \widetilde e^{\vartriangle}(e_1)\neq 0$.
    \end{proof}

    \par The main ingredient of the proof of Theorem \ref{th:surjective} lies in the following lemma: 
    \begin{lemma} \label{lem:linearly}
        Let $E$ be a $p$-natural quasi-Banach space. Given linearly independent $e_1, \hdots, e_n \in E$, there exists $e^{\vartriangle} \in E^{\vartriangle}$ such that $e^{\vartriangle}(e_1), \hdots, e^{\vartriangle}(e_n)$ are linearly independent in $L_p[0,1]$.
    \end{lemma}
    \begin{proof}
    We will proceed by induction in $n$. For $n=2$, let us suppose that, for every $e^{\vartriangle} \in E^{\vartriangle}$, there exist scalars $\lambda_{e^{\vartriangle}}, \mu_{e^{\vartriangle}}$, not both zero, such that $\lambda_{e^{\vartriangle}} e^{\vartriangle}(e_1)+\mu_{e^{\vartriangle}}e^{\vartriangle}(e_2)=0$. Since $E^\vartriangle$ separates the points in $E$, it suffices to prove that there exist scalars $\alpha$ and $\beta$, not both zero, such that $\alpha e^{\vartriangle}(e_1)+\beta e^{\vartriangle}(e_2)=0$ for \emph{every} $e^\vartriangle \in E^{\vartriangle}$.
        
  \par Consider the set
  \[\mathcal{F}=\lbrace e^{\vartriangle}\in E^{\vartriangle}: e^{\vartriangle}(e_2) \neq 0 \rbrace, \]
  and observe that if $e^{\vartriangle} \in \mathcal{F}$ then $\lambda_{e^{\vartriangle}} \neq 0$. Hence we assume that $\lambda_{e^{\vartriangle}}=1$. Let us show that $\mu_{e^\vartriangle}$ can be taken \emph{independently} of $e^{\vartriangle}$; in other words, that there is some $\mu\in \R$ such that $\mu = \mu_{e^{\vartriangle}}$ for every $e^{\vartriangle} \in \mathcal F$. Given $e_1^{\vartriangle}, e_2^{\vartriangle} \in \mathcal{F}$, if $e_1^{\vartriangle}(e_2)$ and $e_2^{\vartriangle}(e_2)$ are linearly independent in $L_p[0,1]$, then
        \begin{align*}
                 -\mu_{e_1^{\vartriangle}} e_1^{\vartriangle}(e_2)-\mu_{e_2^{\vartriangle}} e_2^{\vartriangle}(e_2)=&(e_1^{\vartriangle}+e_2^{\vartriangle})(e_1)=-\mu_{e_1^{\vartriangle}+e_2^{\vartriangle}}(e_1^{\vartriangle}(e_2)+e_2^{\vartriangle}(e_2)) \\
                \iff & (\mu_{e_1^{\vartriangle}+e_2^{\vartriangle}}-\mu_{e_1^{\vartriangle}})e_1^{\vartriangle}(e_2)+(\mu_{e_1^{\vartriangle}+e_2^{\vartriangle}}-\mu_{e_2^{\vartriangle}})e_2^{\vartriangle}(e_2)=0 \\
                \iff & \mu_{e_1^{\vartriangle}+e_2^{\vartriangle}}=\mu_{e_1^{\vartriangle}}=\mu_{e_2^{\vartriangle}}.
            \end{align*}
        On the other hand, if $e_1^{\vartriangle}(e_2)$ and $e_2^{\vartriangle}(e_2)$ are linearly dependent, we can choose $\phi : L_p[0,1] \to L_p[0,1]$ such that $\phi [e_1^{\vartriangle}(e_2)]$ and $e_1^{\vartriangle}(e_2)$ are linearly independent. This of course, forces $\phi [e_1^{\vartriangle}(e_2)]$ and $e_2 ^{\vartriangle}(e_2)$ to be linearly independent as well. Since $\phi\circ e_1^{\vartriangle}\in \mathcal F$, by the previous argument, we obtain that $\mu_{e_1^{\vartriangle}}=\mu_{\phi e_1^{\vartriangle}}=\mu_{e_2^{\vartriangle}}$, as we wanted. 
        
        \par Analogously, in the subset
        $$\mathcal{G}=\lbrace e^{\vartriangle} \in E^{\vartriangle}: e^{\vartriangle}(e_1)\neq 0 \rbrace$$
        we can assume  $\mu_{e^{\vartriangle}}=1$ for all $e^\vartriangle \in E^\vartriangle$, and therefore it follows that there is a privileged scalar $\lambda$ such that $\lambda = \lambda_{e_1^{\vartriangle}}=\lambda_{e_2^{\vartriangle}}$ for every $e_1^{\vartriangle}, e_2^{\vartriangle} \in \mathcal{G}$. Now, since $\mathcal{F} \cap \mathcal{G}\neq \emptyset$ by Lemma \ref{lem:no0}, and every $e^\vartriangle \in \mathcal F \cap \mathcal G$ satisfies
        \begin{equation} \label{eq:linearly}  e^\vartriangle(e_1) + \mu e^\vartriangle(e_2) = \lambda e^\vartriangle(e_1) + e^\vartriangle(e_2) =0, \end{equation}
        the inexorable conclusion is that  $\lambda\mu=1$, and therefore every $e^{\vartriangle} \in \mathcal{F} \cup \mathcal{G}$ satisfies (\ref{eq:linearly}).  On the other hand, if $e^{\vartriangle} \notin \mathcal{F} \cup \mathcal{G}$, then $e^{\vartriangle}(e_1)=0=e^{\vartriangle}(e_2)$, and so the equality $\alpha e^{\vartriangle}(e_1)+\beta e^{\vartriangle}(e_2)=0$ is trivially satisfied for any scalars $\alpha$ and $\beta$.
        
        \par Now, we  prove the general case by a contrapositive argument. Given $n\in \mathbb N$, suppose that $e_1, \ldots, e_{n+1}$ is a collection of linearly independent vectors of $E$ such that $e^{\vartriangle}(e_1), \hdots, e^{\vartriangle}(e_{n+1})$ are linearly dependent for every $e^{\vartriangle} \in E^{\vartriangle}$.

        \begin{claim} If $e^{\vartriangle}(e_i)=0$ for $i\in\{1, \hdots, n\}$, then $e^{\vartriangle}(e_{n+1})=0$.  
        \end{claim}
        \begin{proof}[Proof of the Claim.] \renewcommand{\qedsymbol}{$\blacksquare$} Suppose that a certain $e^\vartriangle \in E^\vartriangle$ satisfies that $e^\vartriangle(e_i)=0$ for every $i\in \{1, \hdots, n\}$, but $e^\vartriangle(e_{n+1})\neq 0$. By the inductive hypothesis, there is $e_0^\vartriangle\in E^\vartriangle$ such that $e_0^\vartriangle(e_1), \hdots, e_0^\vartriangle(e_n)$ are linearly independent vectors in $L_p[0,1]$. Now, choose $x\in L_p[0,1] \setminus \langle e_0^\vartriangle(e_1), \hdots, e_0^\vartriangle(e_n) \rangle$ and pick $\varphi: L_p[0,1] \to L_p[0,1]$ such that $\varphi[e^\vartriangle(e_{n+1})] = x-e_0^\vartriangle(e_{n+1})$. Then the map $e_1^\vartriangle = \varphi \circ e^\vartriangle + e_0^\vartriangle$ satisfies that  ${e_1}^\vartriangle(e_1), \hdots, {e_1}^\vartriangle(e_{n+1})$ are linearly independent, in contradiction with our assumptions. 
        \end{proof} 

        \begin{claim} If $e^\vartriangle(e_i)=0$ for $i\in\{1,\hdots,n-1\}$, then $e^{\vartriangle}(e_n)$ and $e^{\vartriangle}(e_{n+1})$ are linearly dependent.
        \end{claim}
        \begin{proof}[Proof of the Claim.] \renewcommand{\qedsymbol}{$\blacksquare$} If  $e^{\vartriangle} \in E^{\vartriangle}$ is such that $e^{\vartriangle}(e_i)=0$ for all $i\in\{1, \hdots, n-1\}$, and $\lbrace e^{\vartriangle}(e_n), e^{\vartriangle}(e_{n+1})\rbrace$ are linearly independent, then by the strict transitivity of $L_p[0,1]$ we can find $\varphi:L_p[0,1] \to L_p[0,1]$ such that $\varphi(e^{\vartriangle}(e_i))=0$ for every $i\in\{1, \hdots, n\}$, but $\varphi(e^{\vartriangle}(e_{n+1}))\neq 0$, in contradiction with Claim 3. \end{proof}
        
        Now, by appealing to Claim 4 and using the same argument at the beginning of the lemma, one can infer the existence of two scalars $\alpha$ and $\beta$, not both zero, such that $e^\vartriangle(\alpha e_n + \beta e_{n+1})=0$ for \emph{every} $e^\vartriangle\in E^\vartriangle$ such that $e^\vartriangle(e_i)=0$ for $i\in \{1,\hdots, n-1\}$. This contradicts the inductive hypothesis applied to the collection $\lbrace e_1, \hdots, e_{n-1}, \alpha e_n+\beta e_{n+1}\rbrace$. Indeed, if there were $e^\vartriangle\in E^\vartriangle$ such that $e^\vartriangle(e_1), \hdots, e^\vartriangle(e_{n-1}), e^\vartriangle(\alpha e_n + \beta e_{n+1})$ are linearly independent, consider $\varphi: L_p[0,1] \to L_p[0,1]$ such that $\varphi[e^\vartriangle(e_i)]=0$ for every $i\in \{1, \hdots, n-1\}$ but $\varphi[e^\vartriangle(\alpha e_n + \beta e_{n+1})]\neq 0$, and observe that $e_0^\vartriangle = \varphi \circ e^\vartriangle$ satisfies $e_0^\vartriangle(e_i)=0$ for every $i\in \{1,\hdots, n-1\}$ but $e_0^\vartriangle(\alpha e_{n}+\beta e_{n+1})\neq 0$.
    \end{proof}

    \medskip
    \begin{proof}[{\footnotesize PROOF OF THEOREM \ref{th:surjective}}]
        Let us take $f_1, \hdots, f_n \in L_p[0,1]$. Since $e_1, \hdots, e_n$ are linearly independent, by Lemma \ref{lem:linearly}, there exists $e^{\vartriangle} \in E^{\vartriangle}$ such that $e^{\vartriangle}(e_1), \hdots, e^{\vartriangle}(e_n)$ are also linearly independent. Again, using the fact that $L_p[0,1]$ is strictly transitive, exists $\phi:L_p[0,1] \to L_p[0,1]$ such that $\phi(e^{\vartriangle}(e_i))=f_i$, for each $i=1, \hdots, n$. The composition $\phi \circ e^{\vartriangle} \in E^{\vartriangle}$ does the job. 
    \end{proof} 
    \medskip

    \begin{remark} \label{rem:transitive} The above arguments essentially imply that if $E$ is a quasi-Banach space and $F$ is a strictly transitive quasi-Banach space such that $\mathcal L(E,F)$ separate points of $E$, then for every linearly independent subset $\{e_1, .., e_n\}$ in $E$, the evaluation map 
    \[ \mathcal L(E,F) \longrightarrow F^n, \qquad T\longmapsto (T(e_i))_{i=1}^n\]
    is surjective.
    Since $\mathbb R$ is obviously strictly transitive, we deduce that every quasi-Banach space with separating dual is strictly transitive. 
    \end{remark}

	\subsection{Density of $\FVL[E]$ inside $\FpBLp[E]$.} We now apply our results about $p$-natural spaces to show: 
	\begin{prop}
		For every $p$-natural quasi-Banach space $E$, $\FVL[E]$ embeds as a norm dense subspace of $\FpBLp[E]$.
	\end{prop}
    \begin{proof}
        It suffices to check that the vector sublattice generated by $\{\delta_e: e\in E\}$ satisfies the \emph{universal property} of $\FVL[E]$ indicated in Section \ref{section:representation FpBLp}, that is, let us prove that the following commutative diagram holds:
        \begin{equation}
            \notag
            \begin{tikzcd}
        \text{lat}\{\delta_e: e\in E\} \arrow[rd, "\widehat{T}", dashed] &   \\
        E \arrow[u, "\delta_E"] \arrow[r, "T"]                                 & X
        \end{tikzcd}
        \end{equation}
        Let $f$ be a lattice-linear combination of $\delta_{x_1}, \dots, \delta_{x_n}$ for some $x_1,\dots,x_n \in E$, meaning, $f=F(\delta_{x_1},\dots,\delta_{x_n})$ for some lattice-linear expression $F(t_1,\dots,t_n)$. We define $\widehat{T}(f)=F(Tx_1,\dots,Tx_n)$ in $X$ and now check that $\widehat{T}$ is well defined. Suppose that $f=G(\delta_{y_1},\dots,\delta_{y_m})$ where $G(t_1,\dots,t_m)$ is another lattice-linear expression. Let $\lbrace z_1,\dots,z_k \rbrace$ be a maximal linearly independent subset of $\lbrace x_1,\dots,x_n,y_1,\dots,y_m \rbrace$, hence we can write $F(\delta_{x_1},\dots,\delta_{x_n})=\tilde{F}(\delta_{z_1},\dots,\delta_{z_k})$ and $G(\delta_{y_1},\dots,\delta_{y_m})=\tilde{G}(\delta_{z_1},\dots,\delta_{z_k})$. Since $\{\delta_e: e\in E\}$ is a sublattice of $C_{ph}^b(B_{E^{\vartriangle}},L_p[0,1])$, the lattice operations are pointwise, hence
        $$f(e^{\vartriangle})=\tilde{F}(\delta_{z_1}(e^{\vartriangle}),\dots,\delta_{z_k}(e^{\vartriangle}))=\tilde{F}(e^{\vartriangle}(z_1),\dots,e^{\vartriangle}(z_k))$$ 
        in $L_p[0,1]$ for each $e^{\vartriangle} \in B_{E^{\vartriangle}}$. Similarly, $f(e^{\vartriangle})=\tilde{G}(e^{\vartriangle}(z_1),\dots,e^{\vartriangle}(z_k))$. Since $z_1,\dots,z_k$ are linearly independent, Theorem \ref{th:surjective} asserts that $\tilde{F}(t_1,\dots,t_k)=\tilde{G}(t_1,\dots,t_k)$ for all $t_1,\dots,t_k \in L_p[0,1]$. By lattice-linear function calculus (\cite[1.d]{lin-tza}) we obtain that $\tilde{F}(Tz_1, \dots, Tz_k)=\tilde{G}(Tz_1,\dots,Tz_k)$ in $X$, and by the linearity of $T$ it follows that $F(Tx_1,\dots,Tx_n)=G(Ty_1,\dots,Ty_m)$. Therefore, $\widehat{T}$ is well-defined. The uniqueness of the lattice extension $\widehat{T}$ comes from the fact that $\hat{T}\delta_x=Tx$, for every $x \in E$.
    \end{proof}

	\subsection{The functor $\FpBLp[-]$}
    Let us denote by $\QB^{\textbf{(p)}}$ the category of $p$-natural quasi-Banach spaces and operators, and by $\pBL^{\textbf{(p)}}$ the category of $p$-convex $p$-Banach lattices and (lattice) homomorphisms. The relations between categories are usually established by means of \emph{functors}. In our case, we are interested in the definition of the functor associated to the free object $\FpBLp[E]$, for a given $p$-natural quasi-Banach space $E$. Specifically:
		\begin{align*}
			\FpBLp[-]:\QB^{\textbf{(p)}} & \rightsquigarrow \pBL^{\textbf{(p)}} \\
			E & \longmapsto \FpBLp[E] \\
			\lbrace T: E \to F \rbrace & \longmapsto \FpBLp[T]: =\overline{T}:\FpBLp[E] \to \FpBLp[F].
		\end{align*}
        The formal definition of the morphism $\overline{T}:\FpBLp[E] \to \FpBLp[F]$ follows from the universal property of the space $\FpBLp[F]$, as shown in the following diagram:
	\begin{equation}
		\notag
		\begin{tikzcd}
            {\FpBLp[E]} \arrow[drr, "\overline{T}:=\widehat{\delta_{F}T}", bend left=20] & \\
			E \arrow[u, "\delta_E"', hook] \arrow[r, "T"] & 			F \arrow[r, "\delta_{F}", hook] & 			{\FpBLp[F]}                  
		\end{tikzcd}
	\end{equation}
    In other words $\overline{T}:=\widehat{\delta_FT}$.
	\par Observe, nonetheless, that the action of $\overline{T}$ can be concretized by means of the representation of $\FpBLp[E]$. Indeed, since $\FpBLp[E] \subseteq C^b_{ph}(B_{E^{\vartriangle}}, L_p[0,1])$, we can define the composition operator $C_{T^\vartriangle}:\FpBLp[E] \to C^b_{ph}(B_{F^{\vartriangle}}, L_p[0,1])$ given by
	   \[ C_{T^\vartriangle}f(S)=f(T^{\vartriangle}(S)),\]
	for $f \in \FpBLp[E]$, $S \in B_{F^{\vartriangle}}$, and where $T^{\vartriangle}: F^{\vartriangle} \to E^{\vartriangle}$ is the \emph{$\vartriangle$-adjoint} operator of $T:E \to F$ (that is, $T^{\vartriangle}x^{\vartriangle}=x^{\vartriangle}\circ T$, for $x^{\vartriangle}\in F^{\vartriangle}$). In other words, $C_{T^\vartriangle}f=f \circ T^{\vartriangle}$. Let us check that $\overline{T}=C_{T^\vartriangle}$. Note that for any $e \in E$ and every $S \in B_{F^{\vartriangle}}$ we have
    \begin{equation}
    \notag
       \overline{T}\delta_E (e)(S) =\delta_F(Te)(S)=S(Te)=(T^{\vartriangle}S)(e)=\delta_E(e)(T^{\vartriangle}(S))=C_{T^\vartriangle}\delta_E(e)(S).
    \end{equation}
    Hence, as $\overline{T}$ and $C_{T^\vartriangle}$ coincide on all elements of the form $\delta_E(e)$, for $e\in E$, and both are lattice homomorphisms, by uniqueness they must coincide on the whole $\FpBLp[E]$ (and in particular, both map into $\FpBLp[F]$). Finally, note that for the application $C_{T^{\vartriangle}}$ to be well defined we need that $\norm{T}\leq 1$. Since we work with positively homogeneous functions we can always assume that dividing $T$ by its norm when necessary.

\section{Compatibility and automatic extension properties}
    \label{sec:5} 
    
    The approach in this section is motivated by the fact that $p$-Banach spaces are automatically $r$-Banach for every $0<r<p$ and, similarly, $p$-convex quasi-Banach lattices are also $r$-convex for every $0<r<p$. Hence, if $E$ is a $p$-natural quasi-Banach space and $r<p$, then $E$ is also $r$-natural,
    and both free objects $\FxBL{r}[E]$ and $\FpBLp[E]$ could be considered. It is therefore natural to wonder what is the relation between these. In fact, for a natural quasi-Banach space $E$, one might wonder whether the notion of \emph{free $L$-convex quasi-Banach lattice} generated by $E$ can be properly defined and if so, what is the relation to $\FpBLp[E]$ for appropriate $0<p\leq1$.

    \par We shall address both questions at once, and the answer will provide a relation to the convexity of $\FxBL{r}[E]$ when $E$ is a $p$-natural quasi-Banach space and $p>r$. First, let us provide the pertinent definitions.
    
    \par Given $E$ a natural quasi-Banach space, we denote 
    \begin{equation} \label{eq:p_E} p_E = \sup\{ p\in (0,1) : E \text{ is $p$-natural} \}.\end{equation}
    It is worth noting that $E$ is not necessarily $p_E$-natural: take for instance $E=L_{p,\infty}[0,1]$, which satisfies $p_E=p$ but is not $p$-convex \cite{Popa}. However, it is clear that $\FxBL{r}[E]$ can be certainly considered for $r<p_E$. 

       In our attempt to study the relation between $\FpBLp[E]$ and $\FxBL{r}[E]$ for $r,p<p_E$, we will resort to the classical Maurey-Nikishin factorization theorem in the following version --which follows from combining Proposition 43 and Théorème 46(c) with Théorème 8 in \cite{maurey}--: 
    \begin{thm} \cite{maurey} \label{thm:factorization}
        Fix $E$ a natural quasi-Banach space and let $0<p<r<p_E$. Then there exists a constant $C_{p,r}(E)>0$ such that, for every operator $T:E \to L_p(\mu)$, there is:
        \begin{itemize}
            \item an operator $T_0: E \to L_r(\mu)$ with $\|T_0\|\leq C_{p,r}(E)\cdot \|T\|$, 
            \item a non-negative function $g \in L_s(\mu)$ with $\norm{g}_{L_s(\mu)}=1$ and $\frac{1}{p}=\frac{1}{r}+\frac{1}{s}$,
        \end{itemize} 
        so that, for $\varphi_g(f) = g\cdot f$, the following diagram is commutative.
        \begin{equation}
            \notag
            \begin{tikzcd}
E \arrow[rr, "T"] \arrow[rd, "T_0"'] &                                 & L_p(\mu) \\
                                     & L_r(\mu) \arrow[ru, "\varphi_g"'] &         
\end{tikzcd} 
        \end{equation}
    \end{thm}
    \par We now perform an analysis of the constant $ C_{p,r}(E)$ to see that it remains bounded when $p\to 0^+$. For such a purpose, it will be necessary to gather a few known facts about $q$-stable random variables. Given $0<q<2$, a real random variable $X$ is (normalized) \emph{$q$-stable} if $\mathbb E[e^{itX}] = e^{-|t|^q}$ for every $t\in \R$. It is well-known (see e.g. \cite[Theorem 6.4.17]{ak}) that if $X_q$ is a $q$-stable random variable in a probability space $(\Omega, \Sigma, \mu)$ then $X_q\in L_p(\mu)$ precisely when  $0<p<q$. 
    \par The connection between $q$-stable random variables and the Maurey-Nikishin factorization theorem is that Maurey's proof \cite[Lemme 42, Proposition 43 and Théo\-rème 46(c)]{maurey} asserts that, if $E$ is a quasi-Banach space with Rademacher type $q>r$, then Theorem \ref{thm:factorization} is satisfied with  
    \[C_{p,r}(E) \leq T_q(E) \, \frac{A_{r,q}}{A_{p,q}}, \] 
    where $T_q(E)$ is the type $q$ constant of $E$ and $A_{p,q}$ is the $L_p$-norm of a $q$-stable random variable.

    \par Now, if $E$ is a natural quasi-Banach space, then clearly $E$ has type $q$ for any $q<p_E)$. Therefore, given $0<r<p_E$, if we fix any $q\in (r,p_E)$, then to prove that $C_{p,r}(E)$ is uniformly bounded for $p\in (0,r)$ it suffices to show that $A_{p,q}$ is uniformly bounded for such values of $p$. Indeed, the following lemma takes appropriate care of this. It is likely that some parts are well-known, but we include a proof for the sake of completeness. 

    \begin{lemma} \label{lem:q-stable} Given $0<p<q<2$, we have
        \begin{equation} \label{eq:Apq} A_{p,q} = \left[\frac{2\Gamma(p)\Gamma(1-\frac pq)}{\Gamma(\frac p2)\Gamma(1-\frac p2)}\right]^\frac1p,  \end{equation}
        Therefore, 
        \begin{equation}  \label{eq:limit-Apq} \lim_{p\to 0^+} A_{p,q} = e^{\gamma(\frac1q-1)},\end{equation}
        where $\gamma$ denotes the Euler-Mascheroni constant. 
    \end{lemma}
    \begin{proof} If $X_q$ is a $q$-stable random variable and $0<p<q$, then we can write
    \[ \int_\Omega |X_q(\omega)|^p\, d\omega = \int_\R |x|^p d\mu_q(x), \]
    where $\mu_q$ is the probability distribution of $X$. Using the identity
    \[ |x|^p = I(p)^{-1} \int_0^{+\infty} \frac{1-\cos(tx)}{t^{p+1}}\, dt, \qquad I(p) = \int_0^{+\infty} \frac{1-\cos t}{t^{p+1}}\, dt, \]
   which is valid for $0<p<2$, and then applying Fubini's theorem and the fact that $q$-stable random variables are symmetric, we arrive to
    \begin{align*}
        \int_\R |x|^p d\mu_q(x) &= I(p)^{-1} \int_\R \int_0^{+\infty} \frac{1-\cos(tx)}{t^{p+1}}\, dt \, d\mu_q(x) = \\[2mm]
        & = I(p)^{-1} \int_0^{+\infty} \frac{1}{t^{p+1}}\left[1-\int_\R \Re(e^{itx})\, d\mu_q(x)\right] dt = \\[2mm]
        &=I(p)^{-1} \int_0^{+\infty}\frac{1-e^{-t^q}}{t^{p+1}}\, dt. 
    \end{align*}
    Now, observe that the integral
    \[ J(p,q) = \int_0^{+\infty} \frac{1-e^{-t^q}}{t^{p+1}}, \]
    which is clearly convergent for $0<p<q$, can be easily computed by performing integration by parts and then letting $s=t^q$:
    \begin{align*}
        \int_0^{+\infty} \frac{1-e^{-t^q}}{t^{p+1}} \, dt &= \frac{q}{p} \int_0^{+\infty} t^{q-p-1}e^{-t^q} = \frac{1}{p} \int_0^{\infty} s^{ p/q} e^{-s} ds = \frac{\Gamma(1-\frac pq)}{p}.
    \end{align*}
    
    To compute $I(p)$ we employ two standard properties of the Gamma function, namely:
    \begin{enumerate}
        \item[($\dagger$)] For $p>0$, $x^{-p} = \frac{1}{\Gamma(p)}\int_{0}^{+\infty} t^{p-1} e^{-xt}\,  dt$.
        \item [($\ddagger$)] Given $0<a<1$, $\Gamma(a)\cdot\Gamma(1-a) = \int_0^{+\infty} \frac{t^{a-1}}{1+t} \, dt$. 
    \end{enumerate}
    \vspace*{2mm}
    \noindent Now, using first integration by parts, then ($\dagger$) and Fubini's theorem, and finally ($\ddagger$), we obtain
    \begin{align*}
        I(p) &= \frac{1}{p}\int_0^{+\infty} \frac{\sin t}{t^p} \, dt = \frac{1}{p\Gamma(p)}\int_0^{+\infty} \sin t \int_0^{+\infty} t^{p-1}e^{-xt} \, dt\, dx = 
        \\[2mm] &= \frac{1}{p\Gamma(p)} \int_0^{+\infty} t^{p-1} \int_0^{+\infty} \sin x e^{-xt}\, dx\, dt =\\[2mm] &= \frac{1}{p\Gamma(p)}\int_0^{+\infty} \frac{t^{p-1}}{1+t^2} \, dt = \frac{1}{2p\Gamma(p)} \int_0^{+\infty} \frac{s^{p/2-1}}{1+s}\, ds 
        = \\[2mm] &=\frac{\Gamma(\frac p2)\Gamma(1-\frac p2)}{2p\Gamma(p)}. 
    \end{align*}
    This yields (\ref{eq:Apq}). Finally, to compute the limit in (\ref{eq:limit-Apq}), we need to recall a few more facts about the behaviour of the Gamma function: 
    \begin{itemize}
        \item The Laurent series of $\Gamma(z)$ at $z=0$ is 
        \[ \Gamma(z) = \frac1z - \gamma + O(z)\]
        \item The Taylor expansion of $\Gamma(1+z)$ at $z=0$ is
        \[ \Gamma(1+z) = 1-\gamma z + O(z^2).\]
    \end{itemize}
    Hence, by substituting the corresponding expansions, we have 
    \[ \lim_{p\to 0^+} \left[\frac{2\Gamma(p)}{\Gamma(\frac p2)}\right]^\frac1p = \lim_{p\to 0^+} \left[\frac{2(1-\gamma p)}{2-\gamma p}\right]^\frac1p = e^{-\frac\gamma2}, \]
    while, on the other hand, 
    \[ \lim_{p\to 0^+} \left[\frac{\Gamma(1-\frac p2)}{\Gamma(1-\frac pq)} \right]^\frac1p = \lim_{p\to 0^+} \left[\frac{1-\frac\gamma2p}{1-\frac\gamma qp}\right]^\frac1p = e^{\gamma(\frac1q-\frac12)}.\] 
    Multiplying the previous two equations we arrive at (\ref{eq:limit-Apq}). 
    \end{proof} 
    
    The inexorable consequence is: 

    \begin{thm} \label{thm:auto-ext} Let $E$ be a natural quasi-Banach space. Given $0<r<p_E$, there is a constant $\Lambda_r(E)$ such that, for every $0<p<r$, the spaces $\FxBL{p}[E]$ and $\FxBL{r}[E]$ are lattice $\Lambda_r(E)$-isomorphic. 
    \end{thm}
    \begin{proof}
     Fix $p\in (0,r)$, and apply Theorem \ref{thm:factorization} to obtain the following diagram 
        \begin{equation}
            \notag
            \begin{tikzcd}
          & {\FxBL{r}[E]} \arrow[dd, near start, "\widehat{T_0}", dashed] \arrow[rd, "\widehat{T}:=\varphi_g \widehat{T_0}", dashed] &          \\
E \arrow[rr, near start, "T"] \arrow[rd, "T_0"'] \arrow[ru, "\delta", hook] &                                   & L_p[0,1] \\
  & L_r[0,1] \arrow[ru, "\varphi_g"']                                                              &         
\end{tikzcd}
        \end{equation}
        The extension $\widehat{T_0}$ comes automatically from the universal property of $\FxBL{r}[E]$, and $\|\widehat{T_0}\| \leq C_{p,r}(E)\cdot \|T\|$. Now, $\varphi_g$ is a lattice homomorphism, and therefore  the composition $\widehat{T} = \varphi_g \widehat{T}_0$ is a lattice homomorphism extending $T$ and such that         \[
            \|\widehat{T}\| \leq \norm{\varphi_g} \cdot \|\widehat{T_0}\| \leq \|T_0\| \leq C_{p,r}(E) \cdot \norm{T}. 
        \]
        Finally, consider $C_r = \sup_{p\in (0,r)} C_{p,r}$, which is bounded by virtue of Lemma \ref{lem:q-stable}, and observe that every operator $T:E \to L_p[0,1]$ can be extended to a lattice homomorphism $\widehat{T}:\FxBL{r}[E] \to L_p[0,1]$ with $\|\widehat{T}\| \leq C_r(E)\cdot \norm{T}$. Therefore, by the uniqueness of $\FpBLp[E]$, the quasi-Banach lattices $\FpBLp[E]$ and $\FxBL{r}[E]$ must be isomorphic with a constant depending only on $r$. 
    \end{proof}

    \par It is worth to restate Theorem \ref{thm:auto-ext} in the following way:
    
    \begin{cor} \label{cor:auto-ext} Let $E$ be a natural quasi-Banach space. Given $r<p_E$, there exists a constant $C_r(E)$ such that, regardless of $0<p<p_E$, every operator $T: E \to X$, where $X$ is an $r$-convex quasi-Banach lattice, admits a unique extension to a lattice homomorphism $\widehat{T}: \FpBLp[E] \to X$ with $\|\widehat{T}\| \leq C_r(E) \cdot M^{(r)}(E)\cdot \|T\|$.  
    \end{cor}
    
    In the case of Banach spaces, the previous corollary is particularly interesting: 
    \begin{cor} \label{cor:auto-ext-Banach} Let $E$ be a Banach space with non-trivial type. There is a constant $C>0$ (depending only on $E$) such that every operator $T: E \to X$, where $X$ is a $p$-convex $p$-Banach lattice, admits a unique extension to a lattice homomorphism $\widehat{T}: \FBL[E] \to X$ with $\|\widehat{T}\|\leq C \cdot M^{(p)}(X) \cdot \|T\|$. 
	\end{cor}  

    We now show that the previous corollaries are optimal, in the sense that, even if it were the case that $E$ is a $p_E$-natural quasi-Banach space, one cannot take $r=p_E$ in Corollary \ref{cor:auto-ext}. In particular, we show that Corollary \ref{cor:auto-ext-Banach} does not hold when $E$ is a Banach space with trivial type. 
    \par Let $H: L_1(\mathbb R) \to L_{1,\infty}(\mathbb R)$ be (a multiple of) the Hilbert transform 
\[ Hf(x) = \int_\mathbb R \frac{f(t)}{x-t}\, dt\]
where the integral must be understood as a principal value. It is well-known that $H$ is not bounded as an operator from $L_1(\mathbb R)$ to $L_1(\mathbb R)$. For our purposes here, we only need the classical fact that, given $a, b\in \mathbb R$ with $a<b$, then 
\[ 
H\big[I_{[a,b]}\big](x) = \log\left|\frac{x-a}{x-b}\right|, 
\]
where $I_{[a,b]}$ denotes the characteristic function of the interval $[a,b]$. To show that there is no lattice homomorphism  $\widehat{H}: \FBL[L_1(\mathbb R)] \to L_{1,\infty}(\mathbb R)$ such that $\widehat{H}\delta = H$, we observe that 
\[ 
\left\|\sum_{k=1}^n \Big|\delta_{I_{[\frac{k}{n}, \frac{k+1}{n}]}}\Big| \right\|_{\FBL[L_1(\R)]} \leq \sum_{k=1}^n \|I_{[\frac{k}{n}, \frac{k+1}{n}]}\|_{L_1(\R)} = 1.  
\] 
Hence, it is enough to prove that 
the $L_{1,\infty}$-norm of the functions
\[ F_n(x) = \widehat{H}\left[\sum_{k=1}^n \left|\delta_{I_{[\frac{k}{n}, \frac{k+1}{n}]}}\right|\right](x) = \sum_{k=1}^n \left|\log{\left|\frac{x-\frac{k-1}{n}}{x-\frac{k}{n}}\right|}\right| \]
does not remain bounded when $n\to +\infty$. The following properties of $F_n$ will be of use:

\begin{lemma}
    The following hold: 
    \begin{enumerate}
        \item[(i)] $F_n(x) = F_n(1-x)$ for all $x\in \mathbb \R$.
        \item[(ii)] Given $k\in \{0, \hdots, n\}$, $F_n(x)$ is decreasing in $(\frac kn, \frac{2k-1}{2n}]$ and increasing in $[\frac{2k-1}{2n}, \frac{k+1}{n})$. Hence $F_n$ attains local minima at the points $x_{k,n}=\frac{2k-1}{2n}$, $k\in \{1, \hdots, n\}$. 
        \item[(iii)] $F_n(x_{k,n}) \leq F_n(x_{k+1,n})$ provided $k\leq \frac n2$. 
    \end{enumerate}
\end{lemma}
\begin{proof} All assertions follow at once from the identity
\[ F_n(x) = \begin{cases}
  \log\left|\frac{1-x}{x}\right|,  & \text{ if } x<\frac{1}{2n},  \\[2mm]
  \log\left(\frac{x(1-x)}{(x-\frac kn)^2}\right), & \text{ if } \frac{2k-1}{2n} \leq x < \frac{2k-1}{2n}, \\[2mm]
  \log\left|\frac{x}{1-x}\right|, & \text{ if } x \geq \frac{2n-1}{2n},
\end{cases}\]
which follows just by inspection of $F_n$.     
\end{proof}
Now, the previous lemma implies that, when $n\to +\infty$, 
\[ \min_{x\in [0,1]} F_n(x) = F_n\left(\tfrac{1}{2n}\right) = \log(2n-1) \to +\infty, \]
and therefore 
\begin{align*} \|F_n\|_{L_{1,\infty}} &\geq {F_n\!\left(\tfrac{1}{2n}\right)} \cdot\mu\left\{|F_n|>F_n\left(\tfrac{1}{2n}\right)\right\} \geq \log(2n-1) \longrightarrow +\infty. \end{align*}

 \subsection{The free $L$-convex quasi-Banach lattice}
 Let $X$ be an $L$-convex quasi-Banach lattice. Then, there exists $0<p\leq 1$ such that $X$ is $p$-convex, and hence $r$-convex for every $0<r<p$. Moreover, $M^{(r)}(X)$ is a non-decreasing function of $r$. Motivated by these facts, we consider: 
    \[ M^{(L)}(X) = \inf\{M^{(p)}(X): X \text{ is $p$-convex}\}.\]
Observe that $M^{(L)}(L_p[0,1]) = 1$ for all $p\in (0,1]$. Therefore, Proposition \ref{prop:p-natural} asserts that every natural quasi-Banach space is isometrically isomorphic to a subspace of an $L$-convex quasi-Banach lattice $X$ with $M^{(L)}(X)=1$. This motivates the following definition: 

\begin{defi}
	 	Given a natural quasi-Banach space $E$, the \emph{free L-convex quasi-Banach lattice} generated by $E$ is an $L$-convex quasi-Banach lattice $\FqL[E]$, with $M^{(L)}(\FqL[E])=1$,  together with an isometric linear embedding $\delta:E\to\FqL[E]$, such that for every operator $T:E \to X$, where $X$ is an $L$-convex quasi-Banach lattice, there is a unique lattice homomorphism $\widehat{T}:\FqL[E] \to X$ such that $\widehat{T}\delta=T$ and $\|\widehat{T}\|\leq M^{(L)}(X) \cdot \norm{T}$.
	 \end{defi} 
     
     \begin{thm}
        Let $E$ be a natural quasi-Banach space. Then, $\FqL[E]$ exists and it is lattice isomorphic to $\FxBL{r}[E]$ for any $r<p_E$.
    \end{thm}
    \begin{proof} It suffices to show that, for any $r<p_E$ there is a constant $C_r>0$ such that, given any $L$-convex quasi-Banach space $X$ and any norm-one operator $T:E \to X$, there is a lattice homomorphism $\widehat{T}:\FxBL{r}[E] \to X$ extending $T$ and such that $\|\widehat{T} \| \leq C \cdot M^{(L)}(X)$. Hence assume $X$ is $p$-convex for some $p\in (0,r)$. Given $p'\in (0,p]$, Corollary \ref{cor:auto-ext} informs us that, given $p'\in (0,p)$, there is an extension $\widehat{T}: \FxBL{r}[E] \to X$ with $\|\widehat{T}\|\leq C_r \cdot M^{(p')}(X)$. Now, the key fact is that such $\widehat{T}$ is unique, and therefore, it must satisfy $\|\widehat{T}\|\leq C_r \cdot M^{(p')}(X)$ for all $p'\leq p$. Thus we can take infimum and deduce that $\|\widehat{T}\|\leq C_r \cdot M^{(L)}(X)$. 
    \end{proof}

    \begin{cor} \label{cor:pE} $\FqL[E]$ is $p$-convex for every $p\in (0,p_E)$. 
    \end{cor}

    \par To conclude, if $E$ is a $p_E$-natural quasi-Banach space, the preceding arguments are, in principle, no longer applicable. However, we have: 
    \begin{cor}
    Let $E$ be a $p_E$-natural quasi-Banach space. Then $\FqL[E]$ is $p_E$-convex if, and only if, for every $p<p_E$, $\FxBL{p_E}[E]$ and $\FpBLp[E]$ are lattice isomorphic.
    \end{cor}
    \begin{proof} In this proof, we shall denote $\delta_{p_E}: E \hookrightarrow \FxBL{p_E}[E]$ the canonical embedding, to avoid confusion with $\delta: E \hookrightarrow \FqL[E]$. It is clear that $\delta_{p_E}$ can be extended to a lattice homomorphism  $\widehat{\delta_{p_E}}: \FqL[E] \to \FxBL{p_E}[E]$ (see the first diagram of \ref{eq:2-diagrams}). Now, if $\FqL[E]$ is $p_E$-convex, then the canonical embedding $\delta: E \to \FqL[E]$ can also be extended to a lattice homomorphism $\widehat{\delta}: \FxBL{p_E}[E] \to \FqL[E]$ (second diagram of \ref{eq:2-diagrams}). Finally, since $\delta_{p_E}=\widehat{\delta_{p_E}}\delta$ and $\delta=\widehat{\delta}\delta_{p_E}$ it follows that $\widehat{\delta_{p_E}}\widehat{\delta}=\Id_{\FpBL[E]}$ and $\widehat{\delta}\widehat{\delta_{p_E}}=\Id_{\Fq[E]}$. The converse is trivial.     
        \begin{equation}
    \label{eq:2-diagrams}
    \begin{tikzcd}[column sep=1em]
    {\FqL[E]} \arrow[rd, "\widehat{\delta_{p_E}}"]               &            &  & {\FxBL{p_E}[E]} \arrow[rd, "\widehat{\delta}"]             &            \\
    E \arrow[u, "\delta", hook] \arrow[r, "\delta_{p_E}", hook] & {\FxBL{p_E}[E]} &  & E \arrow[r, "\delta", hook] \arrow[u, "\delta_{p_E}", hook] & {\FqL[E]}
    \end{tikzcd}
    \end{equation}
    \end{proof}

\subsection{Applications to $\FpBL$}
\par The above ideas concerning $\FpBLp$ and $\FqL$ can be naturally adapted to shed some light into free quasi-Banach lattices in the non-$p$-convex setting. 

\begin{defi} The \emph{free quasi-Banach lattice} associated to a quasi-Banach space $E$ is a quasi-Banach lattice $\Fq[E]$, together with an isometric embedding $\delta: E \hookrightarrow \Fq[E]$, so that every operator $T$ from $E$ into a quasi-Banach lattice $X$ admits a unique lattice homomorphism $\widehat{T}: \Fq[E] \to X$ such that $\|\widehat{T}\| = \|T\|$.     
\end{defi}

Let $E$ be a quasi-Banach space, and let us denote 
$$p^E=\sup\lbrace p \in (0,1): E \text{ is a } p \text{-Banach space}\rbrace $$
so that $E$ is an $r$-Banach space for $r < p^E$.
Now consider $\Fq[E]$, which is an $s$-Banach space for $s < p^{\Fq[E]}$. Taking $p=\min\lbrace r,s \rbrace$ and arguing as in the proof of Corollary \ref{cor:pE} we obtain the diagrams 
\begin{equation}
    \notag
    \begin{tikzcd}
{\FpBL[E]} \arrow[rd, "\widehat{\delta}"]               &            &  & {\Fq[E]} \arrow[rd, "\widehat{\delta_p}"]             &            \\
E \arrow[u, "\delta_p", hook] \arrow[r, "\delta", hook] & {\Fq[E]},  &  & E \arrow[r, "\delta_p", hook] \arrow[u, "\delta", hook] & {\FpBL[E]}, 
\end{tikzcd}
\end{equation}
whose existence imply that $\Fq[E]$ and $\FpBL[E]$ are lattice isomorphic. The next result follows from the fact that we can repeat this argument for any $q\leq p$. 
\begin{prop}
    Let $E$ be a quasi-Banach space. 
    \begin{enumerate}
        \item If $\Fq[E]$ exists, then $\Fx{p}[E]$ and $\Fx{r}[E]$ are lattice isomorphic for any $p,r < \min\lbrace p^E, p^{\Fq[E]}\rbrace$.
        \item Conversely, if $\Fx{p}[E]$ and $\Fx{r}[E]$ are lattice isomorphic for any $p,r < p^E$, then $\Fq[E]$ exists. 
    \end{enumerate}
\end{prop}

\section{Projectivity in $p$-Banach lattices}
    \label{sec:6}
	In order to proceed, we introduce some categorical notions. Let us denote: 
	\begin{itemize}
		\item $\pB$ the category of $p$-Banach spaces and operators; 
		\item $\pBL$  the category of $p$-Banach lattices and (lattice) homomorphisms.
        \item $\pBLp$ the category of $p$-convex $p$-Banach lattices and (lattice) homomorphisms.
	\end{itemize}
    
    \begin{defi}
        Let $\lambda\geq 1$. An object $P$ of $\pB$ is said to be $\lambda$-projective if given $\pi:X \to Z$ a quotient in $\pB$, for every operator $T:P \to Z$ there is an operator $\tilde{T}:P \to X$ such that $\pi \tilde{T}=T$ and $\| \tilde T\| \leq \lambda\norm{T}$.
    \end{defi}
    We will say that an object $P$ of $\pB$ is \emph{projective} if it is $\lambda$-projective for some $\lambda\geq 1$, and \emph{$1^+$-projective} objects are $(1+\varepsilon)$-projective for all $\varepsilon>0$. It is a well-known fact that the only projective objects in $\pB$ are the $\ell_p(\Gamma)$-spaces (see \cite{ortinsky}). In fact, they are $1^+$-projective.
    \par An analogous definition can be established for projectivity in the $p$-Banach lattice setting.
    \begin{defi}
        Let $\lambda\geq 1$. An object $P$ of $\pBL$ (resp., $\pBLp$) is said to be $\lambda$-projective if given $\pi:X \to Z$ a lattice quotient in $\pBL$ (resp., $\pBLp$), for every lattice homomorphism $T:P \to Z$ there is a lattice homomorphism $\tilde{T}:P \to X$ such that $\pi \tilde{T}=T$ and $\| \tilde T\| \leq \lambda\norm{T}$.
    \end{defi}
    Again, an object $P$ of $\pBL$ (or $\pBLp$) is \emph{projective} if it is $\lambda$-projective for some $\lambda\geq 1$, and \emph{$1^+$-projective} objects are $(1+\varepsilon)$-projective for all $\varepsilon>0$.
    \begin{prop}\label{prop:quotient}
		Every $p$-Banach lattice $X$ of density $\kappa$ is a lattice quotient of $\FpBL[\ell_p(\Gamma)]$, where $\Gamma$ is a set of size $\kappa$. Furthermore, if $X$ is $p$-convex, then it is a lattice quotient of $\FpBLp[\ell_p(\Gamma)]$. 
	\end{prop}
	\begin{proof}
		Given a $p$-Banach lattice $X$ of density $\kappa$, there exists a quotient operator $\pi: \ell_p(\Gamma) \to X$. The extension $\hat{\pi}: \FpBL[\ell_p(\Gamma)] \to X$, given by the universal property of $\FpBL[\ell_p(\Gamma)]$, is the lattice quotient that we need. The argument assuming $p$-convexity is similar.
	\end{proof}
	\begin{prop} \label{prop:projective}
    The space $\FpBL[\ell_p(\Gamma)]$ is a $1^+$-\emph{projective} object in \emph{$\pBL$}.
	 \end{prop}	
	\begin{proof} 
		Let $\pi: X \to Z$ be a lattice quotient between $p$-Banach lattices, $T:\FpBL[\ell_p(\Gamma)] \to Z$ be a lattice homomorphism and $\varepsilon>0$. Since $\ell_p(\Gamma)$ is a $1^+$-projective object in $\pB$ we can always consider an \emph{operator} $\widetilde{T\delta}$ which is a lifting for $T\delta$:
		\[
		\begin{tikzcd}
			X \arrow[r, "\pi"] & Z                                                                                 \\
			& {\FpBL[\ell_p(\Gamma)]} \arrow[u, "T"']                                             \\
			& \ell_p(\Gamma) \arrow[u, "\delta"', hook] \arrow[luu, "\widetilde{T\delta}", bend left]
		\end{tikzcd}
		\]
		with $\|\widetilde{T\delta}\| \leq (1+\varepsilon)\norm{\delta}\norm T \leq (1+\varepsilon)\norm{T}$. Now, we just appeal to the \emph{universal property} of $\FpBL[\ell_p(\Gamma)]$ to obtain a lattice homomorphism $\widehat{T}:\FpBL[\ell_p(\Gamma)] \to X$ such that the following triangle commutes:
		\[
		\begin{tikzcd}
        {\FpBL[\ell_p(\Gamma)]} \arrow[rd, "\widehat{T}"]       &   \\
        \ell_p(\Gamma) \arrow[u, "\delta", hook] \arrow[r, "\widetilde{T\delta}"] & X
        \end{tikzcd}
		\]
        Moreover, since $\delta \widehat{T}=\widetilde{T\delta}$, it holds that $\| \widehat{T}\|\leq (1+\varepsilon)\norm{T}$.
		Finally, we need to check that $\pi\widehat{T}=T$, but this follows from the fact that both operators agree on the generators of $\FpBL[\ell_p(\Gamma)]$; that is, $\pi \widehat{T}\delta=\pi \widetilde{T\delta}=T\delta$.
	\end{proof} 
    Using both Proposition \ref{prop:quotient} and \ref{prop:projective} it is easy to check that:
    \begin{prop}
	    Let $X$ be a $p$-Banach lattice. Then $X$ is a projective object in \emph{$\pBL$} if, and only if, $X$ is lattice complemented in $\FpBL[\ell_p(\Gamma)]$ for some set $\Gamma$. 
	\end{prop}
    \par An analogous proof to the one we have made in the Proposition \ref{prop:projective} works for the case of $p$-convexity. Therefore, it holds that:

        \begin{prop}
            The space $\FpBLp[\ell_p(\Gamma)]$ is a projective object in \emph{${\pBLp}$}.
        \end{prop}
        \subsection{Lattice-projectivity of the spaces $\ell_p(\Gamma)$ using $p$-convexity}
   Thanks to the functional representation of the space $\FpBLp[E]$ when $E$ is a $p$-natural quasi-Banach space, we can deal with some chain conditions in quasi-Banach lattices for the case $E=\ell_p(\Gamma)$, from where we will deduce the projectivity of the spaces $\ell_p(\Gamma)$ in the $p$-convex lattice setting.
	\begin{defi}
	    We say that a vector lattice $X$ satisfies the \emph{countable chain condition} if every family $\mathcal{F} \subset X^+$ such that $f \wedge g = 0$ for all distinct $f,g\in \mathcal{F}$ is countable.  
	\end{defi}
    
    Our next goal is to prove that the space $C_{ph}^b(B_{\ell_p(\Gamma)^{\vartriangle}}, L_p[0,1])$ satisfies the countable chain condition regardless of the set $\Gamma$. We need a couple of lemmata first.
    \begin{lemma} \label{lemma:ccc}
        $B_{\ell_p(\Gamma)^\vartriangle}$ satisfies the (topological) \emph{countable chain condition}; that is to say, every collection of disjoint open sets in $B_{\ell_p(\Gamma)^\vartriangle}$ is countable.
    \end{lemma}
   \begin{proof} There is a natural homeomorphism between $B_{\ell_p(\Gamma)^\vartriangle}$, endowed with the weak-{$\vartriangle$} topology, and $(B_{L_p[0,1]})^\Gamma$, endowed with the product topology given by the usual quasi-normed topology in each $B_{L_p[0,1]}$. 
   Indeed, given $e^\vartriangle \in B_{\ell_p(\Gamma)^\vartriangle}$ we consider the tuple $(x_{\gamma})_{\gamma \in \Gamma}$ given by $x_{\gamma} = e^\vartriangle(e_{\gamma})$, where $(e_\gamma)_{\gamma \in \Gamma}$ is the collection of unit vectors in $\ell_p(\Gamma)$. On the other hand, we associate to a tuple $(x_{\gamma})_{\gamma \in \Gamma} \in (B_{L_p[0,1]})^\Gamma$ an element $e^\vartriangle \in B_{\ell_p(\Gamma)}^{\vartriangle}$ just defining $e^\vartriangle(e_{\gamma})=x_{\gamma}$ for every $\gamma\in \Gamma$.   
   Now, due to the fact that $B_{L_p[0,1]}$ is separable, $(B_{L_p[0,1]})^\Gamma$ enjoys the countable chain condition \cite[2.7.10.d)]{engelking}, and hence so does $B_{\ell_p(\Gamma)^\vartriangle}$.
   \end{proof} 
    
    \begin{lemma}\label{lem:ccc}
        There exists a countable collection of open balls $\mathcal B$ in $L_p[0,1]$ such that:
        \begin{itemize}
            \item[\emph{i)}] $L_p[0,1]_+ \setminus \{0\}\subseteq \bigcup_{B\in \mathcal B} B$.
            \item[\emph{ii)}] Given any set $S$ and two functions $f,g:S \to L_p[0,1]_+$ such that $f^{-1}(B)=g^{-1}(B)$, for any $B\in \mathcal{B}$, then $f=g$.
            \item[\emph{iii)}] If $x,y \in B$, for some $B \in \mathcal B$, then $x \wedge y \neq 0$.
        \end{itemize}
    \end{lemma}
    \begin{proof} Consider  $(b_k)_{k=1}^{\infty}$ a countable dense subset of $L_p[0,1]_+ \setminus \lbrace 0 \rbrace$. It is clear that the collection
        $$\mathcal B = \left\{ B\left( b_k, \frac{\norm{b_k}}{4^n} \right): k, n\in \mathbb N \right\}$$
    satisfies (\emph{i}). 
    \par Let us check $(ii)$. Consider $f,g:S \to L_p[0,1]_+$ two functions such that $f^{-1}(B) = g^{-1}(B)$ for every $B\in \mathcal B$. For a given $x\in L_p[0,1]_+ \setminus \{0\}$ we can find a countable collection of balls $B_{x,n} \in \mathcal B$ such that $\{x\} = \bigcap_{n=1}^\infty B_{x,n}$. Hence $f^{-1}(\lbrace x \rbrace)=g^{-1}(\lbrace x \rbrace)$ for any $x \in L_p[0,1]_+ \setminus \{0\}$, from where it follows that $f=g$. 
        \par Finally, to prove $(\emph{iii)}$, assume $x,y \in B$ for some $B \in \mathcal B$, which is of the form
        $B=B(b_k, \frac{\|b_k\|}{4^n})$ for some naturals $k$ and $n$. Then, the triangle inequality yields
        \[\|x-y\| \leq \|x-b_k\|+\|y-b_k\| \leq \frac{1}{2}\|b_k\| \leq \frac{2}{3}\|x\| < \|x\|, \]
        which is a sufficient condition for $x$ and $y$ not to be disjoint. Indeed, if $x \wedge y=0$, taking norms in the equality $|x-y|=(x \vee y)-(x\wedge y)$, we obtain that
        \[ \norm{x-y} = \norm{|x-y|}= \norm{x \vee y} \geq \norm{x}, \]
        which contradicts our hypothesis.
    \end{proof}
	\begin{thm}\label{thm:ccc}
    The space $C_{ph}^b(B_{\ell_p(\Gamma)^{\vartriangle}}, L_p[0,1])$ satisfies the countable chain condition. 
    \end{thm}
    \begin{proof}
        Let $I$ be an uncountable set and take $(f_i)_{i \in I}$ a family of positive disjoint elements of $C_{ph}^b(B_{\ell_p(\Gamma)^{\vartriangle}}, L_p[0,1])$. Consider $\mathcal B$ a countable collection of balls satisfying properties (i)--(iii) from the previous Lemma. Therefore, 
        for every $i_1, i_2\in I$ such that $i_1\neq i_2$, there exists an open ball $B_{k(i_1,i_2)}\in\mathcal B$ such that $f^{-1}_{i_1}(B_{k(i_1,i_2)}) \neq f^{-1}_{i_2}(B_{k(i_1,i_2)})$. Therefore, we can find a subfamily of elements $(f_j)_{j \in J}$, with $J \subset I$ and $|J|=|I|$, so that there exists an open ball $B \in \mathcal B$ satisfying that $f^{-1}_{j_1}(B) \neq f^{-1}_{j_2}(B)$, for any $j_1, j_2 \in J$. 
        Now, observe that the open sets $(f^{-1}_j(B))_{j \in J}$ are disjoint, thanks to (iii) of Lemma \ref{lem:ccc}. But, by virtue of Lemma \ref{lem:ccc}, $B_{\ell_p(\Gamma)^{\vartriangle}}$ enjoys the countable chain condition, and so this implies that $J$ --and therefore, $I$-- is a countable set. 
    \end{proof}
	
	\begin{thm} The space $\ell_p(\Gamma)$ is a projective object in \emph{$\pBL$} if, and only if, the set $\Gamma$ is countable. 
	\end{thm}
	\begin{proof}
		First, let us show that $\ell_p(\Gamma)$ is projective whenever $\Gamma$ is countable. Since every complemented sublattice of a projective object in $\pBL$ is also a projective object in $\pBL$, it suffices to show that $\ell_p$ can be embedded as a lattice-complemented subspace of $\FpBL[\ell_p]$. First, observe that the identity map in $\ell_p$ induces a lattice projection $\beta: \FpBL[\ell_p] \to \ell_p$. We now define a homomorphism $\alpha: \ell_p \to \FpBL[\ell_p]$ such that $\beta \alpha = \operatorname{Id}_{\ell_p}$.   
		Let us consider the map
		\[ \alpha: \ell_p \to \FpBL[\ell_p], \quad \alpha(e_n) = \left[|\delta_{e_n}|-4^n\left(\sum_{j<n}|\delta_{e_j}|+\sum_{j>n} 2^{-j}|\delta_{e_j}|\right)\right]_+
		\]
		which clearly satisfies $\beta \alpha = \operatorname{Id}_{\ell_p}$. To check that $\alpha$ is a norm-one lattice homomorphism with closed range, denote $\alpha(e_n)=f_n$, and observe that, since $\|f_n\|\leq 1$, one has:
		\[
		\norm{\sum_{n=1}^{\infty}a_n\alpha(e_n)}=\norm{\sum_{n=1}^{\infty}a_nf_n} \leq \left(\sum_{n=1}^{\infty} |a_n|^p \right)^{1/p}.
		\]
		On the other hand, since $\beta(f_n)=e_n$, we obtain that
		\[
		\left( \sum_{n=1}^{\infty} |a_n|^p \right)^{1/p}=\norm{\sum_{n=1}^{\infty}a_ne_n}_{\ell_p}=\norm{\sum_{n=1}^{\infty} a_n \beta(f_n)} \leq \norm{\sum_{n=1}^{\infty} a_nf_n}=\norm{\sum_{n=1}^{\infty}a_n\alpha(e_n)}.\]
		This is enough to show that $\ell_p$ is a projective object in $\pBL$.
		\par Conversely, assume that $\Gamma$ is uncountable. To prove that $\ell_p(\Gamma)$ is not projective we will make use of $\FpBLp[\ell_p(\Gamma)]$. Consider a lattice quotient $\pi: \FpBLp[\ell_p(\Gamma)] \to \ell_p(\Gamma)$. If $\ell_p(\Gamma)$ were projective, we could consider a lifting of the identity map of $\ell_p(\Gamma)$ as in the following diagram:
		\begin{equation}
			\notag
			\begin{tikzcd}
				{\FpBLp[\ell_p(\Gamma)]} \arrow[r, "\pi"] & \ell_p(\Gamma)                                                           \\
				& \ell_p(\Gamma) \arrow[u, equal] \arrow[lu, "\widetilde{\Id}"]
			\end{tikzcd}            
		\end{equation}
		Let us observe that the map $\widetilde{\Id}$ is necessarily injective, and therefore, $(\widetilde{\Id}(e_i))_{i\in I}$ is a collection of pairwise disjoint vectors in $\FpBLp[\ell_p(\Gamma)]$, which contradicts Theorem \ref{thm:ccc}. Therefore, $\ell_p(\Gamma)$ is not a projective object in $\pBLp$ when $\Gamma$ is uncountable, so neither can it be in $\pBL$.
	\end{proof}	
    \begin{cor}
        The space $\ell_p(\Gamma)$ is a projective object in \emph{$\pBLp$} if, and only if, $\Gamma$ is countable.
    \end{cor}

\section{Miscellaneous open questions and remarks} 

We close the paper with some problems related to the material presented above which we have not been able to solve: 

\subsection{A functional representation for $\FpBL[E]$.} Perhaps the most annoying question we have not been able to solve is to provide a functional representation for $\FpBL[E]$ for any $p$-Banach space $E$. The main obstacle we found is that it is not clear that there exists a privileged $p$-Banach lattice $X_p$ such that $\mathcal L(E,X_p)$ separates the points in $E$ regardless of the $p$-Banach space $E$. 
Indeed, in \cite[Theorem 10.3]{kalton-linear}, N. Kalton shows that (a real version of) the space $L_p(\mathbb T)/H_p$ does not admit non-zero operators to any $\sigma$-complete order continuous quasi-Banach lattice.

\subsection{$\FpBLp[E]$ and the countable chain condition}
We conjecture that Theorem \ref{thm:ccc} is true for \emph{any} $p$-natural quasi-Banach space $E$; that is, $C_{ph}^b(B_{E^\vartriangle}, L_p[0,1])$ satisfies the countable chain condition regardless of $E$. The arguments of Section \ref{sec:6} can be carefully adapted to show that this is certainly the case when $E$ is separable. Indeed, since there is a quotient map $\pi: \ell_p \to E$, then the map $\pi^\vartriangle: {E^\triangle} \hookrightarrow \ell_p^\vartriangle$ is $\sigma(E^\vartriangle,E)$-$\sigma(\ell_p^\vartriangle, \ell_p)$-continuous. Hence, $B_{E^\triangle}$ is a subspace of the separable and metrizable space $B_{\ell_p^\triangle} = (B_{L_p[0,1]})^\mathbb N$ and therefore, $B_{E^\triangle}$ has the (topological) countable chain condition, which is enough to conclude. However, since the topological countable chain condition is not hereditary, this argument does not seem to admit a straigthforward generalization.

\subsection{Maurey-Nikishin factorization revisited}
In Section \ref{sec:5} we showed that $\FqL[L_1(\mathbb R)]$ is not lattice isomorphic to $\FBL[L_1(\mathbb R)]$. It is clear that the construction presented there is also valid for $L_1[0,1]$. However, we have no candidate for an operator $T: L_p[0,1] \to L_{p,\infty}[0,1]$ which cannot be extended to a lattice homomorphism from $\FpBLp[L_p[0,1]]$ to $L_{p,\infty}[0,1]$. 
\par Observe that, by Theorem \ref{thm:auto-ext}, every such operator does admit an extension $\widehat{T}: \FxBL{r}[L_p[0,1]] \to L_{p,\infty}[0,1]$ provided $0<r<p$. On the other hand, the Maurey-Nikishin factorization theorem through $L_{p,\infty}(\mu)$ --see \cite[Theorem III.H.6.]{woj}-- hints at the possibility that, if every operator $T:L_p(\mu) \to L_{p,\infty}(\mu)$ can be extended to a lattice homomorphism $\widehat{T}: \FpBLp[L_p(\mu)] \to L_{p,\infty}(\mu)$, then so can every operator $T: L_p(\mu) \to L_0(\mu)$. 

\par In general, one can wonder whether there is a quasi-Banach space which is $p_E$-natural and such that $\FxBL{p_E}[E]$ is not lattice isomorphic to $\FxBL{r}[E]$ for any $r<p_E$. Indeed, the existence of and operator $T: L_p[0,1] \to L_{p,\infty}[0,1]$ which cannot be extended to a lattice homomorphism from $\FpBLp[L_p[0,1]]$ would imply that $E=L_p[0,1]$ affirmatively solves this question.

\section*{Acknowledgements}
We wish to thank J. L. Ansorena for bringing reference \cite{turpin} to our attention. This meeting happened in the context of the conference ``Structures in Banach spaces'', which took place at the Erwin Schrödinger Institute (Vienna) during March 2025. We thank the organizers for making it possible. 
\par 
This research has been supported by grants PID2020-116398GB-I00, PID2024-162214NB-I00 and CEX2023-001347-S funded by  MCIN/AEI/10.13039/501100011033. The first author has also been supported by the grant PID2023-146505NB-C21 funded by MICIU/AEI/
10.13039/501100011033. The third author benefited from a predoctoral grant associated to ``Programa de Financiación UCM-Banco Santander (CT24/25)'', awarded by Universidad Complutense de Madrid and Banco Santander.

    \bibliography{bibliography}{}
\bibliographystyle{plain}

\end{document}